\documentclass[a4paper]{amsart}
\usepackage{amssymb} 
\usepackage{amsfonts}
\usepackage{amsxtra}
\usepackage[T1]{fontenc}
\usepackage[latin1]{inputenc}
\usepackage{ae}
\usepackage[all]{xy}
\usepackage{url}
\usepackage{verbatim}
\usepackage{pgf}
\usepackage{tikz}
\usepackage{tikz-cd}
\usepackage{stmaryrd}

\usetikzlibrary{arrows,patterns}

\include{diagram}

\newcommand*{\ket}{\rangle}
\newcommand*{\bra}{\langle}

\newcommand*{\BA}{{\bf A}}
\newcommand*{\B}{\mathcal{B}}
\newcommand*{\C}{\mathcal{C}}

\newcommand*{\D}{\mathcal{D}}
\newcommand*{\E}{\mathcal{E}}
\newcommand*{\F}{\mathcal{F}}

\renewcommand*{\O}{\mathcal{O}}
\renewcommand*{\P}{\mathcal{P}}

\newcommand*{\BB}{{\bf B}}
\newcommand*{\BU}{{\bf U}}

\newcommand*{\BV}{{\bf V}}

\newcommand*{\BW}{{\bf W}}

\newcommand*{\BX}{{\bf X}}
\newcommand*{\BY}{{\bf Y}}
\newcommand*{\BZ}{{\bf Z}}
\newcommand*{\Bb}{{\bf b}}
\newcommand*{\Bc}{{\bf c}}
\newcommand*{\Bf}{{\bf f}}
\newcommand*{\Bg}{{\bf g}}
\newcommand*{\Bh}{{\bf h}}
\newcommand*{\Bi}{{\bf i}}

\newcommand*{\Bk}{{\bf k}}

\newcommand*{\HH}{\mathbb{H}}
\newcommand*{\KH}{\mathbb{K}}
\newcommand*{\LH}{\mathbb{L}}
\newcommand*{\Split}{\mathfrak{S}}
\newcommand*{\Rep}{\mathsf{Rep}}

\newcommand*{\Hilb}{\mathsf{Hilb}}
\newcommand*{\HILB}{\mathsf{HILB}}

\newcommand*{\Ob}{\mathsf{Ob}}

\DeclareMathOperator{\Bicat}{Bicat}

\DeclareMathOperator{\Cat}{Cat}

\DeclareMathOperator{\im}{im}
\DeclareMathOperator{\ind}{ind}
\DeclareMathOperator{\Ind}{Ind}

\DeclareMathOperator{\id}{id}

\DeclareMathOperator{\Lin}{Lin}
\DeclareMathOperator{\lin}{lin}
\DeclareMathOperator{\Bilin}{Bilin}

\newenvironment{bnum}
{\begin{list}{}
    {\setlength{\labelwidth}{15pt}
     \setlength{\leftmargin}{\labelwidth}
    }
}
{\end{list}}

\numberwithin{equation}{section}
\theoremstyle{change}
\newtheorem{theorem}{Theorem}[section]
\newtheorem{prop}[theorem]{Proposition}
\newtheorem{lemma}[theorem]{Lemma}

\newtheorem{definition}[theorem]{Definition}
\newtheorem{example}[theorem]{Example}

\begin{document}

\title{On bicolimits of $ C^* $-categories}

\author{Jamie Antoun}
\address{School of Mathematics and Statistics \\
         University of Glasgow \\
         University Place \\
         Glasgow G12 8SQ\\
         United Kingdom 
}
\email{j.antoun.1@research.gla.ac.uk}

\author{Christian Voigt}
\address{School of Mathematics and Statistics \\
         University of Glasgow \\
         University Place \\
         Glasgow G12 8SQ\\
         United Kingdom 
}
\email{christian.voigt@glasgow.ac.uk}



\maketitle

\begin{abstract}
We discuss a number of general constructions concerning additive $ C^* $-categories, focussing in particular on establishing the existence of bicolimits. As 
an illustration of our results we show that balanced tensor products of module categories over $ C^* $-tensor categories exist without any finiteness assumptions. 
\end{abstract}

\section{Introduction} 

The study of $ C^* $-tensor categories and their module categories has intimate connections with quantum groups, subfactors, and quantum field theory, see for 
instance \cite{DRduality}, \cite{Wenzlcstartensorqg}, \cite{DecommerYamashitaTK1}, \cite{PVcstartensor}, \cite{BKLRbook}. It can be viewed as an incarnation 
of \emph{categorified algebra} in a framework adapted to operator algebras. In this context it is useful to be able to perform a number of constructions 
with $ C^* $-categories, possibly equipped with further structure, in analogy to $ C^* $-algebras and their representations. 

Basic examples of such constructions arise from tensor products. In the purely algebraic setting, one works with the Deligne tensor product of $ k $-linear 
abelian categories \cite{Delignetannaka}, or more generally, the Kelly tensor product of finitely cocomplete $ k $-linear categories over a field $ k $, 
see \cite{LopezFrancotensorproducts}, \cite{DSPSbalancedtensorproduct}. Balanced tensor products of module categories appear for instance in the study of categorical 
Morita equivalence \cite{Muegersft1}, \cite{ENOfusionhomotopy} and in topological field theory \cite{BZBJsurfaces}. While many applications are concerned with 
semisimple module categories over fusion categories, it is known that balanced tensor products, and more generally bicolimits, exist in much greater 
generality. This follows from abstract results in enriched category theory, using that the categories under consideration are categories of algebras over a 
finitary $ 2 $-monad \cite{BKPmonad}, and includes weighted bicolimits. 
It is not obvious, however, to what extent this machinery can be adapted to the setting of $ C^* $-categories, since some of the ingredients used 
in \cite{BKPmonad} are no longer available in this case. 

Strict colimits of $ C^* $-categories were studied by Dell'Ambrogio \cite{DellAmbrogiounitarymodel}, who gave a construction based on generators 
and relations. More recently, an approach towards the construction of bicolimits of $ C^* $-categories has been taken by Albandik 
and Meyer, using the language of $ C^* $-correspondences. In their paper \cite{AlbandikMeyercolimits} they discuss several concrete examples of bicolimits 
in this setting, and prove an abstract existence result under additional assumptions.  
Compared to the more general theory of $ C^* $-categories, an advantage of the bicategory of $ C^* $-correspondences is that all ingredients are very concrete, 
and that correspondences are well-adapted for applications to Cuntz-Pimsner algebras and generalised crossed products \cite{AlbandikMeyerore}. However, 
from the point of view of higher category theory it is more natural to work at the level of $ C^* $-categories, which allows one at the same time to bypass 
some limitations of the techniques used in \cite{AlbandikMeyercolimits}. 

The main aim of this paper is to show that certain $ 2 $-categories of additive $ C^* $-categories are closed under conical bicolimits. 
Here the term \emph{additive} means that we require our $ C^* $-categories to admit certain types of direct sums. 
For the sake of definiteness we will mainly focus our attention on the $ 2 $-category $ C^* \Lin $ of countably additive $ C^* $-categories, 
that is, $ C^* $-categories admitting countable direct sums. Compared to the widely used notion of a finite direct sum, infinite direct 
sums are less familiar in the context of $ C^* $-categories. In fact, their behaviour differs significantly from the purely algebraic situation. 

It should be noted that $ C^* $-categories are not particularly well-behaved with respect to abstract categorical concepts in general. For instance, in the 
absence of strong finiteness assumptions they rarely admit kernels or cokernels. One peculiarity is that the $ * $-structure allows one to reverse the direction 
of arrows, so that existence of a certain type of limits is equivalent to existence of the corresponding type of colimits. As a consequence, a number of 
standard concepts need to be modified in order to remain meaningful. When it comes to infinite direct sums this leads one almost inevitably to work with 
non-unital $ C^* $-categories. Accordingly, multiplier categories appear from the very start. 

Our argument for establishing the existence of bicolimits is a variant of the adjoint functor theorem, albeit in the setting of bicategories instead of 
ordinary categories. Compared to the considerations in \cite{AlbandikMeyercolimits}, this approach has the advantage of being 
applicable in quite general circumstances. One could adapt it to give direct existence proofs for bicolimits in other types 
of $ 2 $-categories, like the $ 2 $-category of finitely cocomplete $ k $-linear categories over a field $ k $. 
However, the construction only gives limited information on how bicolimits look concretely. 

As an illustration we discuss the construction of balanced tensor products of module categories over $ C^* $-tensor categories,  
without having to impose semisimplicity or rigidity assumptions, or other finiteness conditions. 
While not surprising, already this special case of our general existence result seems not at all obvious 
from the outset, given the nuances in the $ C^* $-setting mentioned above. 

Let us now explain how the paper is organised. In section \ref{secpreliminaries} we have collected various definitions and constructions related 
to $ C^* $-categories. This includes in particular a careful discussion of direct sums, and of ind-categories 
of $ C^* $-categories. We also introduce a notion of finitely presentable objects in a $ C^* $-category, and of finitely accessible $ C^* $-categories. 
Hilbert modules provide a rich source of examples, and we review the link between the correspondence bicategory and singly 
generated $ C^* $-categories. In section \ref{sectensorproducts} we recall the construction of minimal and maximal tensor product of $ C^* $-categories. 
We discuss their functoriality and provide a new characterisation of the maximal tensor norm. 
The main result of the paper is contained in section \ref{secbicolim}, where we prove the existence of bicolimits in the $ 2 $-category $ C^* \Lin $ of 
countably additive $ C^* $-categories. Finally, in section \ref{secbalancedtensorproducts} we study $ C^* $-tensor categories and their module 
categories in $ C^* \Lin $. Based on our main result we show that balanced tensor products always exist in this setting, and give some explicit examples. 
For the convenience of the reader we have assembled a few definitions and facts regarding bicategories in the appendix.  

Let us conclude with some remarks on our notation and conventions. 
Unless explicitly stated otherwise we assume all categories to be small. 
If $ \BV $ is a category we write $ \BV(V,W) $ for the set of morphisms between objects $ V, W \in \BV $. The space of adjointable operators between 
Hilbert modules $ \E, \F $ is denoted by $ \LH(\E, \F) $, and we write $ \KH(\E, \F) $ for the subspace of compact operators. 
The closed linear span of a subset $ X $ of a Banach space is denoted by $ [X] $. Depending on the context, 
the symbol $ \otimes $ denotes the algebraic tensor product over the complex numbers or various completions thereof. We sometimes 
write $ \odot $ for algebraic tensor products.  
By slight abuse of language, all our $ C^* $-categories are semi-categories by default, that is, they do not necessarily contain identity morphisms. 
We will speak of unital categories or $ 1 $-categories if we are dealing with categories in the usual sense.

\section{Preliminaries} \label{secpreliminaries}

In this section we review definitions and results regarding $ C^* $-categories and fix our notation. 

\subsection{$ C^* $-categories} 

We begin with the definition of $ C^* $-categories and their basic properties. For additional information see for 
instance \cite{GLRwstar}, \cite{Mitchenercstarcategories}. 

By a $ * $-category we shall mean a semicategory $ \BV $ such that all morphism spaces $ \BV(V,W) $ for $ V, W \in \BW $ are complex vector spaces and 
the composition maps $ \BV(X,Y) \times \BV(Y, Z) \rightarrow \BV(X, Z), (f,g) \mapsto g \circ f $ are bilinear, together with an antilinear involutive 
contravariant endofunctor $ *: \BV \rightarrow \BV $ which is the identity on objects, mapping $ f \in \BV(X, Y) $ to $ f^* \in \BV(Y, X) $. 

A $ C^* $-category is a $ * $-category $ \BV $ such that all morphism spaces $ \BV(V,W) $ are complex Banach spaces, the composition 
maps $ \BV(X,Y) \times \BV(Y, Z), (f,g) \mapsto g \circ f $ satisfy $ \|g \circ f\| \leq \|g\| \|f \| $, the $ C^* $-identity 
$$
||f^* \circ f || = ||f ||^2 
$$ 
holds and $ f^* \circ f \in \BV(X,X) $ is positive for all $ f \in \BV(X, Y) $. 
One may phrase this as saying that the category $ \BV $ is enriched in the category of complex Banach spaces and contractive linear maps 
with some further structure and properties. By a unital $ C^* $-category we mean a $ C^* $-category which contains identity morphisms for all its objects. 

A basic example of a (large) $ C^* $-category is the category $ \HILB = \HILB_\mathbb{C} $ of Hilbert spaces with morphisms all \emph{compact} linear operators 
between them. More generally, for any $ C^* $-algebra $ A $ we have the large $ C^* $-category $ \HILB_A $ of right Hilbert $ A $-modules with compact 
adjointable operators as morphisms. 

A $ * $-functor $ \Bf: \BV \rightarrow \BW $ between $ * $-categories is a functor such that $ \Bf(f^*) = \Bf(f)^* $ for all morphisms $ f $. 
A linear $ * $-functor is a $ * $-functor such that all associated maps on morphism spaces are linear. If $ \Bf: \BV \rightarrow \BW $ is a 
linear $ * $-functor between $ C^* $-categories then the maps $ \Bf: \BV(V, W) \rightarrow \BW(\Bf(V), \Bf(W)) $ are automatically contractive. 
A natural transformation $ \tau: \Bf \Rightarrow \Bg $ of linear $ * $-functors $ \Bf, \Bg: \BV \rightarrow \BW $ between $ C^* $-categories is a uniformly 
bounded family of morphisms $ \tau(V): \Bf(V) \rightarrow \Bg(V) $ such that $ \tau(W) \circ \Bf(f) = \Bg(f) \circ \tau(V) $ for all morphisms $ f \in \BV(V,W) $. 
Linear $ * $-functors between $ C^* $-categories together with their natural transformations as morphisms form naturally $ C^* $-categories.  
A $ * $-linear functor $ \Bf: \BV \rightarrow \BW $ between unital $ C^* $-categories is called unital if $ \Bf(\id_X) = \id_{\Bf(X)} $ for all $ X \in \BV $. 
A natural transformation between unital $ * $-linear functors is defined in the same way as in the non-unital case. 

Let $ X, Y \in \BV $ be objects in a $ C^* $-category. A left multiplier morphism $ L: X \rightarrow Y $ is a collection of uniformly bounded 
linear maps $ L(Z): \BV(Z, X) \rightarrow \BV(Z,Y) $ such that $ L(W)(h \circ g) = L(Z)(h) \circ g $ for all $ h \in \BV(Z, X) $ and $ g \in \BV(W, Z) $. 
Similarly, a right multiplier morphism $ R: X \rightarrow Y $ is a collection of uniformly bounded linear maps $ R(Z): \BV(Y, Z) \rightarrow \BV(X, Z) $ 
such that $ R(Z)(f \circ h) = f \circ R(W)(h) $ for all $ f \in \BV(W, Z) $ and $ h \in \BV(Y, W) $. A \emph{multiplier morphism} $ M: X \rightarrow Y $ 
is a pair $ (L, R) $ of left and right multiplier morphisms from $ X $ to $ Y $ such that $ g \circ L(W)(f) = R(Z)(g) \circ f $ for all $ f \in \BV(W, X) $ 
and $ g \in \BV(Y, Z) $. Clearly, every morphism $ f: X \rightarrow Y $ in $ \BV $ defines a multiplier morphism $ M_f = (L_f, R_f): X \rightarrow Y $ by 
setting $ L_f(h) = f \circ h, R_f(g) = g \circ f $ for $ h \in \BV(W, X) $ and $ g \in \BV(Y, Z) $. 

If $ F = (L_F, R_F), G = (L_G, R_G) $ are multiplier morphisms from $ X $ to $ Y $ and $ Y $ to $ Z $, respectively, then the 
composition $ G \circ F = (L_G \circ L_F, R_F \circ R_G) $ is a multiplier morphism from $ X $ to $ Z $. 
If $ F = (L_F, R_F): X \rightarrow Y $ is a multiplier morphism then the adjoint $ F^*: Y \rightarrow X $ is the multiplier 
morphism $ F^* = (L_{F^*}, R_{F^*}) $ where $ L_{F^*}(f) = R_F(f^*)^* $ and $ R_{F^*}(f) = L_F(f^*)^* $. The identity maps define naturally a multiplier 
morphism $ \id_X: X \rightarrow X $ for $ X \in \BV $. 

If we write $ M \BV(X,Y) $ for the set of all multiplier morphisms from $ X $ to $ Y $,  
then we obtain a category $ M\BV $, the multiplier category of $ \BV $, with the same objects as $ \BV $ and morphism sets $ M \BV(X,Y) $. The multiplier 
category $ M\BV $ is naturally enriched over Banach spaces such that $ F = (L_F, R_F) \in M\BV(X,Y) $ has norm $ \|F\| = \sup(\|L_F\|, \|R_F\|) $ where 
$$
\|L_F\| = \sup_{Z \in \BV} \|L_F(Z)\|, \qquad \|R_F\| = \sup_{Z \in \BV} \|R_F(Z)\|. 
$$ 
With this structure the category $ M\BV $ becomes a unital $ C^* $-category, compare \cite{Kandelakimultiplier}, \cite{Vassellibundles}. 
Two objects $ X, Y $ in a $ C^* $-category $ \BV $ are said to be isomorphic if they are isomorphic in $ M\BV $ in the usual sense. 

The multiplier category $ M\HILB $ has as objects all Hilbert spaces and morphisms the bounded linear operators between them. Similarly, 
the unital $ C^* $-category $ M\HILB_A $ for a $ C^* $-algebra $ A $ has the same objects as $ \HILB_A $ but all adjointable linear operators as morphisms. 

Let $ \BV $ be a $ C^* $-category and $ \BX, \BY \in \BV $. The \emph{strict topology} on $ M\BV(X,Y) $ is the locally convex topology defined by saying that 
$ g_i \rightarrow g $ strictly iff $ g_i \circ f \rightarrow g \circ f $ and $ h \circ g_i \rightarrow h \circ g $ in norm for all $ f \in \BV(W,X), h \in \BV(Y,Z) $. 
The space of multiplier morphisms $ M\BV(X,Y) $ can be viewed as the completion of $ \BV(X,Y) $ with respect to the strict topology. 
A linear $ * $-functor $ \Bf: \BV \rightarrow M\BW $ is called strict if the associated maps $ \BV(X,Y) \rightarrow M\BW(\Bf(X),\Bf(Y)) $ are all strictly 
continuous on bounded subsets. This is equivalent to saying that for every $ X \in \BV $ and approximate identity $ (e_i)_{i \in I} $ in $ \BV(X,X) $ 
the net $ \Bf(e_i)_{i \in I} $ converges strictly in $ M\BW(\Bf(X), \Bf(X)) $. 
Every strict linear $ * $-functor $ \Bf: \BV \rightarrow M\BW $ extends uniquely to a strict linear $ * $-functor $ M\BV \rightarrow M\BW $,  
which we will again denote by $ \Bf $. 
If $ \BV $ is unital then every linear $ * $-functor $ \Bf: \BV \rightarrow M\BW $ is strict. 

If $ \BV $ is a $ * $-category and $ \BW $ a $ C^* $-category then a linear $ * $-functor $ \Bf: \BV \rightarrow M\BW $ is called nondegenerate if 
\begin{align*}
[\Bf(\BV(Y,Y)) \circ \BW(\Bf(X), \Bf(Y))] &= \BW(\Bf(X), \Bf(Y)) = [\BW(\Bf(X), \Bf(Y)) \circ \Bf(\BV(X,X))] 
\end{align*}
for all $ X, Y \in \BV $. 
If $ \BV $ is a $ C^* $-category then a nondegenerate linear $ * $-functor $ \Bf: \BV \rightarrow M\BW $ is strict and extends uniquely to a unital 
linear $ * $-functor $ M\BV \rightarrow M\BW $. 
The composition of nondegenerate linear $ * $-functors between $ C^* $-categories is again nondegenerate. 

Given nondegenerate linear $ * $-functors $ \Bf, \Bg: \BV \rightarrow M\BW $, a multiplier natural transformation $ \tau: \Bf \Rightarrow \Bg $ is a family of 
uniformly bounded multiplier morphisms $ \tau(V) \in M\BW(\Bf(V), \Bg(V)) $ satisfying the usual commutativity with respect to morphisms $ f \in \BV(V,W) $. 
If $ \BV $ is a $ C^* $-category then we automatically get $ \tau(W) \circ \Bf(f) = \Bg(f) \circ \tau(V) $ for all $ f \in M\BV(V,W) $, using the unique unital 
extensions $ \Bf, \Bg: M\BV \rightarrow M\BW $. 
A (unitary) natural isomorphism $ \tau: \Bf \Rightarrow \Bg $ between nondegenerate linear $ * $-functors is a multiplier natural transformation such 
that $ \tau(V) $ is a (unitary) isomorphism for all $ V \in \BV $. We write $ \Bf \cong \Bg $ if there exists a natural isomorphism from $ \Bf $ to $ \Bg $. 
In this case there exists also a unitary natural isomorphism $ \Bf \Rightarrow \Bg $, see Proposition 2.6 in \cite{DellAmbrogiounitarymodel}. 

Two $ C^* $-categories $ \BV, \BW $ are called \emph{equivalent} if there exist nondegenerate linear $ * $-functors $ \Bf: \BV \rightarrow M\BW $ 
and $ \Bg: \BW \rightarrow M\BV $ such that their mutual compositions are naturally isomorphic to the identities. 
In this case we actually have $ \Bf: \BV \rightarrow \BW, \Bg: \BW \rightarrow \BV $, and we write $ \BV \simeq \BW $. 

A linear $ * $-functor $ \Bf: \BV \rightarrow M\BW $ is called faithful if the associated maps $ \Bf: \BV(X,Y) \rightarrow M\BW(\Bf(X), \Bf(Y)) $ are 
injective for all $ X, Y \in \BV $. By a \emph{realisation} of a $ C^* $-category $ \BV $ we shall mean a nondegenerate faithful 
linear $ * $-functor $ \pi: \BV \rightarrow M\HILB $ which is injective on objects. 
Every $ C^* $-category admits a realisation, compare Proposition 1.14 in \cite{GLRwstar} and Theorem 6.12 in \cite{Mitchenercstarcategories}.

\subsection{Direct sums and subobjects} 

Let us discuss the notion of a direct sum in a $ C^* $-category. 

\begin{definition} \label{defdirectsum}
Let $ \BV $ be a $ C^* $-category and let $ I $ be a set. If $ (V_i)_{i \in I} $ is a family of objects in $ \BV $, then a direct sum of $ (V_i)_{i \in I} $ 
is an object $ \bigoplus_{i \in I} V_i \in \BV $ together with multiplier morphisms $ \iota_j: V_j \rightarrow \bigoplus_{i \in I} V_i $ for $ j \in I $ 
such that 
$$ 
\iota_k^* \circ \iota_l = \delta_{kl} \id_{V_l}, \qquad \sum_{j \in I} \iota_j \circ \iota_j^* = \id, 
$$ 
where the sum in the second expression is required to converge in the strict topology. 
\end{definition} 

We may also express the conditions in Definition \ref{defdirectsum} in terms of $ \iota_j^* = \pi_j: \bigoplus_{i \in I} V_i \rightarrow V_j $ in $ M\BV $. 
If the index set $ I = \{1, \dots, n\} $ is finite we also write 
$$ 
\bigoplus_{i = 1}^n V_i = V_1 \oplus \cdots \oplus V_n 
$$ 
for a direct sum. Note that the multiplier morphisms $ \iota_j: V_j \rightarrow \bigoplus_{i = 1}^n V_i $ 
satisfy $ \iota_1 \circ \iota_1^* + \cdots + \iota_n \circ \iota_n^* = \id $ in this case, and that we have a canonical identification 
$$
M\BV(\bigoplus_{i = 1}^n V_i, \bigoplus_{i = 1}^n V_i) = \bigoplus_{i,j = 1}^n M\BV(V_i, V_j).  
$$
In particular, if $ \BV $ is a unital $ C^* $-category then finite direct sums are automatically both products and coproducts in the sense of category theory. 
By definition, a direct sum indexed by the empty set is a zero object, that is an object $ 0 \in \BV $ such that $ \BV(0,V) = 0 = \BV(V, 0) $ 
for all $ V \in \BV $. 

Let $ I $ be again arbitrary. If $ f_i: V_i \rightarrow W $ is a uniformly bounded family of morphisms such that $ \sum_{i \in I} f_i \circ \pi_i $ converges 
strictly then there exists a unique multiplier morphism $ f: \bigoplus_{i \in I} V_i \rightarrow W $ such that $ f \circ \iota_i = f_i $ for all $ i \in I $. 
Using this mapping property one checks that a direct sum of $ (V_i)_{i \in I} $ is unique up to isomorphism.  
Note that a nondegenerate linear $ * $-functor $ \Bf: \BV \rightarrow M\BW $ preserves all direct sums which exist in $ \BV $. 
 
The $ C^* $-category $ \HILB_A $ of Hilbert modules over a $ C^* $-algebra $ A $ admits arbitrary direct sums, 
given by direct sums of Hilbert $ A $-modules in the usual sense. That is, we have 
$$
\bigoplus_{i \in I} V_i = \{(v_i)_{i \in I} \in \prod_{i \in I} V_i \mid \sum_{i \in I} \bra v_i, v_i \ket \text{ converges in } A \}
$$
for a family $ (V_i)_{i \in I} $ of Hilbert $ A $-modules. 

If $ \BV $ is a $ C^* $-category one can form a $ C^* $-category $ \BV^{\oplus} $, the \emph{finite additive completion} 
of $ \BV $, by formally adjoining to $ \BV $ all finite direct sums of objects in $ \BV $. More explicitly, the objects 
in $ \BV^{\oplus} $ are families $ (V_i)_{i \in I} $ of objects $ V_i \in \BV $ indexed by some finite set $ I $.  
To construct the morphism spaces in $ \BV^{\oplus} $ we use a realisation $ \iota: \BV \rightarrow M\HILB $, and 
let $ \BV^{\oplus}((V_i)_{i \in I}, (W_j)_{j \in J}) $ be the closed subspace 
$$
\bigoplus_{i \in I, j \in J} \iota(\BV(V_i, W_j)) 
\subset \LH(\bigoplus_{i \in I} \iota(V_i), \bigoplus_{j \in J} \iota(W_j)) 
$$
with the induced algebraic operations. There is an obvious nondegenerate linear $ * $-functor $ \BV \rightarrow \BV^{\oplus} $ obtained by viewing objects 
of $ \BV $ as indexed by a one-element index set. The finite additive completion is closed under taking finite direct sums. 
In a similar way one can define additive completions with respect to arbitrary regular cardinals. 

Let $ \BV $ be a $ C^* $-category and $ V \in \BV $. A \emph{subobject} of $ V $ is an object $ U \in \BV $ together with a multiplier 
morphism $ \iota: U \rightarrow V $ such that $ \iota^* \circ \iota = \id_U $. In this case $ p = \iota \circ \iota^* \in M\BV(V,V) $ is a projection.
Conversely, we say that a projection $ p \in M\BV(V,V) $ is split if there exists a subobject $ U $ of $ V $ with $ \iota: U \rightarrow V $ 
such that $ p = \iota^* \circ \iota $. The $ C^* $-category category $ \BV $ is called \emph{subobject complete} if every projection in $ \BV $ is split. 

The \emph{subobject completion} $ \Split(\BV) $ of a $ C^* $-category $ \BV $ is the category whose objects are pairs $ (X, p) $ where $ X \in \BV $ 
and $ p \in M\BV(X,X) $ is a projection. Morphisms from $ (X, p) $ to $ (Y, q) $ in $ \Split(\BV) $ are all morphisms $ h \in \BV(X, Y) $ 
satisfying $ h \circ p = h = q \circ h $. Composition of morphisms and the $ * $-structure are inherited from $ \BV $. 
It is straightforward to check that $ \Split(\BV) $ is a $ C^* $-category in a natural way, and that $ \Split(\BV) $ is subobject complete. 
If $ \BV $ is closed under taking direct sums of some cardinality, then the same is true for $ \Split(\BV) $.

\subsection{Additive $ C^* $-categories} 

We will be interested in $ C^* $-categories which admit at least finite direct sums, and mainly focus on $ C^* $-categories which are closed under 
countable direct sums and subobject complete. 
For convenience we shall introduce the following terminology.  

\begin{definition} 
A countably additive $ C^* $-category is a subobject complete small $ C^* $-category which admits all countable direct sums. 
By a finitely additive  $ C^* $-category we mean a subobject complete unital small $ C^* $-category which admits all finite direct sums. 
\end{definition} 

We write $ C^* \Lin $ for the $ 2 $-category which has countably additive $ C^* $-categories as objects, all nondegenerate 
linear $ * $-functors $ \Bf: \BV \rightarrow M\BW $ as $ 1 $-morphisms from $ \BV $ to $ \BW $, together with their multiplier natural transformations 
as $ 2 $-morphisms. Similarly, we let $ C^* \lin $ be the $ 2 $-category which has finitely additive $ C^* $-categories as objects, unital 
linear $ * $-functors $ \Bf: \BV \rightarrow \BW $ as $ 1 $-morphisms, and their natural transformations as $ 2 $-morphisms. 

Note that every object $ \BV $ in either of these $ 2 $-categories has a zero object $ 0 $, and that $ \BV(V,W) $ contains at least the zero morphism 
for any $ V, W \in \BV $. 
Let us also point out that there are no nontrivial unital categories in $ C^* \Lin $ since the existence of a nonzero object and infinite direct sums always 
gives rise to objects with nonunital endomorphism algebras. 
By the same token, categories in $ C^* \lin $ are not closed under taking infinite direct sums, with the trivial exception being the zero category, just 
containing a zero object. 

A standard example of a category in $ C^* \lin $ is the category $ \Hilb^f_A $ of finitely generated projective Hilbert $ A $-modules over a 
unital $ C^* $-algebra $ A $. 
The objects of $ \Hilb^f_A $ are Hilbert modules isomorphic to direct summands in $ A^{\oplus n} $ for some $ n \in \mathbb{N} $, and the morphisms are all 
adjointable operators between them. Note that all adjointable operators between finitely generated projective Hilbert $ A $-modules are automatically compact. 
Similarly, a prototypical example of a category in $ C^* \Lin $ is the category $ \Hilb_A $ of Hilbert modules  
over an arbitrary $ C^* $-algebra $ A $ which are isomorphic to direct summands of the standard Hilbert module $ \HH_A = \bigoplus_{n = 1}^\infty A $, 
with all compact linear operators as morphisms. In either case we tacitly need to make arrangements to ensure that the categories under considerations are small 
by choosing a set of such Hilbert modules which is large enough to accommodate the constructions we want to consider. 

By construction both $ \Hilb_A $ and $ \Hilb_A^f $ are full subcategories of $ \HILB_A $. If $ A $ is $ \sigma $-unital, then $ \Hilb_A $ consists precisely of 
all countably generated Hilbert $ A $-modules by Kasparov's stabilisation theorem. 
We note that Kasparov's stabilisation theorem may fail for Hilbert modules which are not countably generated. This is related to the (non-) existence of frames, 
see \cite{Linoframes}. 

A \emph{generator} for $ \BV \in C^* \Lin $ is an object $ G \in \BV $ such that any object $ V \in \BV $ is isomorphic to a direct summand 
of $ \bigoplus_{n = 1}^\infty G $. We say that $ \BV $ is singly generated if it admits a generator. Clearly the $ C^* $-category $ \Hilb_A $ is singly 
generated by $ A $, viewed as Hilbert module over itself. 

\begin{prop} \label{singlygenerated}
Let $ \BV \in C^* \Lin $ be singly generated. Then $ \BV \simeq \Hilb_A $ for some $ C^* $-algebra $ A $. 
\end{prop} 

\proof Let $ G $ be a generator for $ \BV $ and write $ A = \BV(G,G) $. We construct a linear $ * $-functor $ F: \BV \rightarrow M\Hilb_A $ by 
setting $ F(V) = \BV(G, V) $, viewed as Hilbert $ A $-module with inner product $ \bra f, g \ket = f^* \circ g $ and the module structure given by 
right multiplication. The functor $ F $ is easily seen to be nondegenerate, 
so that $ F $ commutes with direct sums.  
Restricted to the full subcategory of $ \BV $ formed by at most countable direct sums of copies of $ G $, the functor $ F $ factorises through $ \Hilb_A $ and 
is fully faithful by construction. 
Every object of $ \BV $ is isomorphic to a direct summand of $ \bigoplus_{n \in \mathbb{N}} G $, 
which implies that $ F: \BV \rightarrow \Hilb_A $ is fully faithful. Since $ \BV $ is closed under subobjects the functor $ F $ is essentially surjective. \qed 

Recall that an $ A $-$ B $-correspondence for $ C^* $-algebras $ A, B $ is a Hilbert $ B $-module $ \E $ together with a 
nondegenerate $ * $-homomorphism $ \phi: A \rightarrow \LH(\E) $, compare \cite{BMZhighercategory}. As in \cite{AlbandikMeyercolimits} 
we shall say that $ \E $ is proper if the image of $ \phi: A \rightarrow \LH(\E) $ is contained in $ \KH(\E) $. 
By definition, an $ A $-$ B $-correspondence in $ \Hilb_B $ is an $ A $-$ B $-correspondence whose underlying Hilbert $ B $-module is contained in $ \Hilb_B $. 

Every $ A $-$ B $-correspondence $ \E $ in $ \Hilb_B $ defines a nondegenerate linear $ * $-functor $ - \otimes_A \E: \Hilb_A \rightarrow M\Hilb_B $, 
sending $ V \in \Hilb_A $ to the interior tensor product $ V \otimes_A \E $. 
If $ \E $ is proper $ A $-$ B $-correspondence in $ \Hilb_B $ then $ - \otimes_A \E $ defines a nondegenerate linear $ * $-functor $ \Hilb_A \rightarrow \Hilb_B $. 
In particular, every $ * $-homomorphism $ f: A \rightarrow B $ between separable $ C^* $-algebras $ A $ and $ B $ induces a 
nondegenerate linear $ * $-functor $ \Hilb_A \rightarrow \Hilb_B $ by considering the Hilbert $ B $-module $ \E = [f(A)B] $ with the left action of $ A $ 
induced by $ f $. Note that $ f $ is not required to be nondegenerate. 

Conversely, we have the following result. 

\begin{prop} \label{hilbmodfunctors}
Let $ A, B $ be $ C^* $-algebras. Then every nondegenerate linear $ * $-functor $ \Bf: \Hilb_A \rightarrow M\Hilb_B $ 
is of the form $ \Bf \cong - \otimes_A \E $ for some $ A $-$ B $-correspondence $ \E \in \Hilb_B $. 
Similarly, every nondegenerate linear $ * $-functor $ \Bf: \Hilb_A \rightarrow \Hilb_B $ 
is of the form $ \Bf \cong - \otimes_A \E $ for a proper $ A $-$ B $-correspondence $ \E \in \Hilb_B $. 
\end{prop} 

\proof Let $ \Bf: \Hilb_A \rightarrow M\Hilb_B $ be a nondegenerate linear $ * $-functor. The Hilbert $ B $-module $ \E = \Bf(A) \in \Hilb_B $ is 
equipped with a nondegenerate left action of $ A $ via $ A \cong \KH(A) \rightarrow \LH(\Bf(A)) = \LH(\E) $. Let $ \Bg: \Hilb_A \rightarrow M\Hilb_B $ 
be the nondegenerate linear $ * $-functor given by $ \Bg(V) = V \otimes_A \E $. We have a canonical unitary 
isomorphism $ \sigma(A): \Bg(A) = A \otimes_A \E \rightarrow \E = \Bf(A) $. 
Both $ \Bf $ and $ \Bg $ are nondegenerate and hence compatible with direct sums, so that we obtain a canonical unitary 
isomorphism $ \sigma(V): \Bg(V) \rightarrow \Bf(V) $ for $ V = \HH_A $ as well. 
Since every object of $ \Hilb_A $ is isomorphic to a subobject of $ \HH_A $, it follows that $ \sigma $ extends uniquely to a unitary natural 
isomorphism $ \Bf \Rightarrow \Bg $. 

If $ \Bf $ is a nondegenerate linear $ * $-functor from $ \Hilb_A $ into $ \Hilb_B $, then the $ A $-$ B $ correspondence $ \E = \Bf(A) \in \Hilb_B $ is 
clearly proper. \qed 

Proposition \ref{hilbmodfunctors} shows in particular that the bicategory of $ C^* $-algebras and $ C^* $-correspondences \cite{BMZhighercategory} 
fits naturally into the framework of countably additive $ C^* $-categories. Apart from the fact that we allow more $ 2 $-morphisms, a technical 
difference is that we obtain size restrictions on the correspondences appearing in our setup. This point is however of little practical relevance. 

Note that $ \Hilb_A $ determines $ A $ only up to Morita equivalence, so that working at the level of $ C^* $-categories provides a ``coordinate free'' 
description of the Morita class of the $ C^* $-algebra $ A $ in the correspondence bicategory. 
Using Proposition \ref{singlygenerated} and Proposition \ref{hilbmodfunctors} one can of course translate between 
the correspondence bicategory and singly generated countably additive $ C^* $-categories.

\subsection{Ind-categories of $ C^* $-categories} 

In this paragraph we discuss how the definition of ind-categories can be adapted from the purely algebraic setting in order to be compatible 
with $ * $-structures. This extends the considerations for semisimple unital $ C^* $-categories in \cite{NYdrinfeldcenter}. 

Recall that a partially ordered set $ I $ is called $ \kappa $-directed for some regular cardinal $ \kappa $ if all subsets of $ I $ of cardinality less 
than $ \kappa $ admit an upper bound. Our constructions below could be done at this level of generality, but we shall restrict ourselves to directed sets 
in the usual sense, which are precisely the $ \kappa $-directed sets for $ \kappa = \aleph_0 $. 

Let $ \BV $ be a $ C^* $-category and let $ I $ be a directed set. By an inductive system over $ I $ in $ \BV $, or ind-object, we mean an inductive
system $ X = ((X_i)_{i \in I}, (\iota_{ij})_{i \leq j \in I}) $ in the multiplier category $ M\BV $ such that all connecting 
morphisms $ \iota_{ij}: X_i \rightarrow X_j $ are isometries. That is, the morphisms $ \iota_{ij} \in M\BV $ are required to satisfy $ \iota_{ii} = \id $ 
and $ \iota_{jk} \circ \iota_{ij} = \iota_{ik} $ for all $ i \leq j \leq k $. 
We will often abbreviate this as $ X = (X_i)_{i \in I} $ with the connecting morphisms suppressed.  

Let $ X = ((X_i)_{i \in I}, (\iota_{ik})_{i \leq k \in I}) $ and $ Y = ((Y_j)_{j \in J}, (\eta_{jl})_{j \leq l \in J}) $ be ind-objects in $ \BV $. We define 
$$
\Ind \BV(X, Y) = \varinjlim_i \varinjlim_j \BV(X_i, Y_j) 
$$
to be the Banach space inductive limit of the inductive system $ \BV(X_i, Y_j) $ 
with respect to the isometric connecting maps $ \BV(X_i, Y_j) \rightarrow \BV(X_k, Y_l) $ for $ i \leq k, j \leq l $ 
given by $ f \mapsto \eta_{jl} \circ f \circ \iota_{ik}^* $. This differs from the definition of ind-morphisms in \cite{NYdrinfeldcenter}, and the appearance 
of inductive limits in both variables may appear strange at first sight. We will explain further below how this relates to the approach in \cite{NYdrinfeldcenter}. 

If $ Z = ((Z_k)_{k \in K}, (\kappa_{km})_{k \leq m \in K}) $ is another ind-object then the composition maps 
$ \BV(X_i, Y_j) \times \BV(Y_j, Z_k) \rightarrow \BV(X_i, Z_k) $ in $ \BV $ induce a well-defined composition 
$ \Ind \BV(X, Y) \times \Ind \BV(Y, Z) \rightarrow \Ind \BV(X, Z) $. To this end it suffices to observe 
that elements $ f \in \Ind \BV(X, Y), g \in \Ind \BV(Y, Z) $ in the algebraic inductive limits may be represented 
by morphisms $ f_{ij} \in \BV(X_i, Y_j), g_{jk} \in \BV(Y_j, Z_k) $ by choosing $ i,j,k $ large enough, and that the class of the composition 
$ g_{jk} \circ f_{ij} $ does not depend on the choice of $ j $. 
There is also a canonical $ * $-operation on morphism spaces sending $ f \in \BV(X_i, Y_j) \subset \Ind \BV(X, Y) $ 
to $ f^* \in \BV(Y_j, X_i) \subset \Ind \BV(Y, X) $. 

\begin{definition} 
Let $ \BV $ be a $ C^* $-category. The ind-category $ \Ind \BV $ of $ \BV $ is the subobject completion of the $ C^* $-category with objects all 
ind-objects over $ \BV $ and morphism spaces $ \Ind \BV(X, Y) $ as defined above. We write $ \ind \BV $ for the full subcategory of $ \Ind \BV $ consisting 
of all subobjects of countable inductive systems in $ \BV $.  
\end{definition} 

The $ C^* $-category $ \BV $ embeds into $ \Ind \BV $ by considering objects of $ \BV $ as constant ind-objects indexed by a one-element set, and the embedding 
functor $ \BV \rightarrow \Ind \BV $ is fully faithful on morphisms. 

Let us discuss the structure of morphism spaces in the multiplier category of $ \Ind \BV $. Consider 
objects $ X = ((X_i)_{i \in I}, (\iota_{ik})_{i \leq k \in I}) $ and $ Y = ((Y_j)_{j \in J}, (\eta_{jl})_{j \leq l \in J}) $ in $ \Ind \BV $, and let us 
define the space of \emph{formal multiplier morphisms} from $ X $ to $ Y $ by 
\begin{align*}
\lim_\leftrightarrows M\BV(X, Y) 
&= \{(f_{ij}) \mid \sup_{i,j} \|f_{ij}\| < \infty, \eta_{jk}^* \circ f_{ik} = f_{ij} \forall j \leq k, f_{lj} \circ \iota_{il} = f_{ij} \forall i \leq l \} \\
&\subset \prod_{(i,j) \in I \times J}  M\BV(X_i, Y_j). 
\end{align*}
We say that $ f = (f_{ij}) \in \lim_\leftrightarrows M\BV(X, Y) $ is \emph{right strict} iff for all $ i \in I, g \in \BV(X_i, X_i) $ and $ \epsilon > 0 $ 
there exists $ j_0 \in J $ such that 
$$
\|\eta_{jl} \circ f_{ij} \circ g - f_{il} \circ g\| < \epsilon 
$$
for all $ l \geq j \geq j_0 $. 
Note that this is automatically satisfied if $ Y $ is a constant inductive system. 
Similarly, let us say that $ f = (f_{ij}) \in \lim_\leftrightarrows M\BV(X, Y) $ is \emph{left strict} iff for all $ j \in J, h \in \BV(Y_j, Y_j) $ 
and $ \epsilon > 0 $ there exists $ i_0 \in I $ such that 
$$
\|h \circ f_{ij} \circ \iota_{ik}^* - h \circ f_{kj} \| < \epsilon 
$$
for all $ k \geq i \geq i_0 $. This is automatic if $ X $ is constant. 

We define 
\begin{align*}
M\text{-}\varprojlim_i \varinjlim_j M\BV(X_i, Y_j) = \{(f_{ij})_{i \in I, j \in j} \mid f \text{ left strict and right strict} \} 
\end{align*} 
as a subspace of $ \lim_\leftrightarrows M\BV(X, Y) $. If $ X \in \BV $ is viewed as a constant ind-object we will also use the notation 
\begin{align*}
M\text{-}\varinjlim_j M\BV(X, Y_j) 
= \{(f_j)_{j \in j} \mid f \text{ right strict} \} 
\end{align*}
for this space and call it the \emph{multiplier inductive limit}. 
In a dual fashion, we write $ M\text{-}\varprojlim_j M\BV(Y_j, X) $ for the \emph{multiplier projective limit}, 
which is obtained by taking pointwise adjoints in $ M\text{-}\varinjlim_j M\BV(X, Y_j) $. Equivalently, 
\begin{align*}
M\text{-}\varprojlim_j M\BV(Y_j, X) 
= \{(f_j)_{j \in j} \mid f \text{ left strict} \} 
\end{align*} 
as a subspace of the Banach space projective limit $ \varprojlim_j M\BV(Y_j, X) $. 

\begin{prop} 
Let $ \BV $ be a $ C^* $-category. Then there exists a canonical isometric linear isomorphism 
$$ 
M \Ind \BV(X, Y) \rightarrow M\text{-}\varprojlim_i \varinjlim_j M\BV(X_i, Y_j) 
$$ 
for all $ X = (X_i)_{i \in I}, Y = (Y_j)_{j \in J} \in \Ind \BV $. 
\end{prop} 

\proof Note first that if $ Z = (Z_k)_{k \in K} $ is an ind-object and $ l \in K $ then the identity in $ M\BV(Z_l, Z_l) $ induces 
a canonical multiplier morphism $ \iota_l: Z_l \rightarrow Z $ in $ M \Ind \BV(Z_l, Z) $. Here we consider $ Z_l $ as a constant ind-object. 

In order to prove the assertion let us start with the case that $ X $ is a constant inductive system and $ Y = (Y_j)_{j \in J} $ arbitrary. 
For $ f \in M \Ind \BV(X, Y) $ we get multiplier morphisms $ f_j = \iota_j^* \circ f \in M\BV(X, Y_j) $ for all $ j \in J $. 
The family $ (f_j)_{j \in J} $ is contained in $ M\text{-} \varinjlim_j M\BV(X, Y_j) $, 
and we obtain a contractive linear map $ \phi: M \Ind \BV(X, Y) \rightarrow M\text{-} \varinjlim_j M\BV(X, Y_j) $ in this way. Assume that $ \phi(f) = 0 $ and 
let $ h \in \BV(X,X) $. Then $ \iota_j^* \circ f \circ h = 0 $, and so $ \iota_j \circ \iota_j^* \circ f \circ h = 0 $ for all $ j \in J $. 
On the other hand we have $ \iota_j \circ \iota_j^* \circ f \circ h \rightarrow f \circ h $ since $ f \circ h \in \Ind \BV(X, Y) $. 
Hence $ f \circ h = 0 $, and it follows that $ f = 0 $. This shows that $ \phi $ is injective. 
Next assume $ (f_j) \in M\text{-}\varinjlim_j M\BV(X, Y_j) $ and let $ h \in \BV(W, X) $.
Factorising $ h = g \circ k $ with $ g \in \BV(X,X) $ we see that for any $ \epsilon > 0 $ there exists $ j_0 \in J $ such 
that $ \|\eta_{jl} \circ f_j \circ h - f_l \circ h \| < \epsilon $ for $ l \geq j \geq j_0 $. 
Hence we obtain a uniquely determined element $ f \circ h \in \Ind \BV(W, Y) $ with $ \lim_j \iota_j \circ f_j \circ h = f \circ h $. 
If $ Z \in \BV $ we also define $ k \circ f $ for $ k \in \varinjlim_j \BV(Y_j, Z) = \Ind \BV(Y, Z) $ by $ k \circ f = \lim_j k \circ \iota_j \circ f_j $. 
This is well-defined since 
\begin{align*}
\|k \circ \iota_j \circ f_j - k \circ \iota_l \circ f_l \| 
&= \|k \circ \iota_j \circ \eta_{jl}^* \circ f_l  - k \circ \iota_l \circ f_l \| \\
&\leq \|k \circ \iota_j \circ \eta_{jl}^* - k \circ \iota_l \| \|f\| < \epsilon 
\end{align*}
for all $ l \geq j $ sufficiently large. 
Both these constructions extend to morphisms in $ \Ind \BV $, 
and one checks that $ f $ defines a multiplier morphism in $ M \Ind \BV(X, Y) $ such that $ \phi(f) = (f_j) $. 
Hence $ \im(\phi) = M\text{-}\varinjlim_j M\BV(X, Y_j) $, and the map $ \phi: M \Ind \BV(X, Y) \rightarrow M\text{-} \varinjlim_j M\BV(X, Y_j) $ 
is in fact an isometric isomorphism. 

Now consider $ f \in M \Ind \BV(X,Y) $ for arbitrary $ X, Y \in \Ind \BV $. We get morphisms $ f_{ij} = \iota_j^* \circ f \circ \iota_i \in M\BV(X_i, Y_j) $ for 
all $ i \in I, j \in J $, and  by our above considerations the family $ (\iota_j^* \circ f \circ \iota_i)_{j \in J} $ defines an element $ f \circ \iota_i $ 
in $ M \text{-} \varinjlim M\BV(X_i, Y) $ for each $ i \in I $. Since $ (X_i)_{i \in I} $ is an inductive system these morphisms assemble to an 
element of the Banach space projective limit $ \varprojlim_i M \text{-} \varinjlim_j M\BV(X_i, Y_j) $, and 
the resulting linear map $ \phi: M \Ind \BV(X, Y) \rightarrow \varprojlim_i M \text{-} \varinjlim_j M\BV(X_i, Y_j) $ is injective. 
Moreover the image of $ \phi $ is contained in $ M \text{-} \varprojlim_i \varinjlim_j M\BV(X_i, Y_j) $, since for $ f \in M \Ind \BV(X, Y) $ 
and $ h \in \BV(Y_j, Y_j) $ we have $ \|h \circ \iota_j^* \circ f \circ \iota_i \circ \eta^*_{ik} - h \circ \iota_j^* \circ f \circ \iota_k \| < \epsilon $ 
for all $ k \geq i \geq i_0 $ with sufficiently large $ i_0 $, using that $ h \circ \iota_j^* \circ f \in \Ind \BV(X, Y_j) $. 

Assume $ (f_{ij}) \in M \text{-} \varprojlim_i \varinjlim_j M\BV(X_i, Y_j) $. Due to our previous considerations, the $ i $-th component of this family
gives a multiplier morphism $ f_i \in M \Ind \BV(X_i, Y) $ for every $ i \in I $. Moreover, if $ W \in \BV $ and $ g \in \Ind \BV(W, X) $ 
then there exists a unique morphism $ f \circ g \in \Ind \BV(W, Y) $ such that $ \lim_i f_i \circ \iota_i^* \circ g = f \circ g $. 
This assignment extends to define a multiplier morphism $ f \in M \Ind \BV(X,Y) $. We obtain $ \phi(f) = (f_{ij}) $, and the map $ \phi $ is surjective and isometric. 
\qed 

Assume that $ \BV \in C^* \lin $ is a $ C^* $-category such that all morphism spaces in $ \BV $ are finite dimensional. Then $ M\BV(X_i, Y_j) = \BV(X_i, Y_j) $ 
for all $ i,j $, and the multiplier category $ M \Ind \BV $ 
agrees with the ind-category defined in this setting in \cite{NYdrinfeldcenter}. More precisely, Lemma 2.1 in \cite{NYdrinfeldcenter} shows that in this 
case multiplier projective and inductive limits reduce to ordinary Banach space projective and inductive limits, respectively. 

\begin{definition} \label{definductivelimit}
Let $ \BV $ be a $ C^* $-category. The direct limit of an inductive system $ (X_i)_{i \in I} $ in $ \BV $ is an object $ \varinjlim_i X_i $ in $ \BV $ 
together with a family of isometries $ \iota_j \in M\BV(X_j, \varinjlim_i X_i) $ for $ j \in I $ such that $ \iota_j \circ \iota_{ij} = \iota_i $ for 
all $ i \leq j $ and $ \iota_j \circ \iota_j^* \rightarrow \id $ in the strict topology. 
\end{definition} 

Let us emphasize that a direct limit in the sense of Definition \ref{definductivelimit} is typically not a direct limit in the categorical sense. 
Since we will not deal with categorical direct limits this should not lead to confusion. 

If $ X = \varinjlim_{i \in I} X_i $ is a direct limit then $ \BV(X,X) $ is the closed linear span of all morphisms of the form $ \iota_j \circ f \circ \iota_i^* $ 
for $ f \in \BV(X_i, X_j) $ and $ i,j \in I $. 
In general, if $ (f_i) $ is a family of multiplier morphisms in $ M\BV(X_i, Y) $ satisfying $ f_j \circ \iota_{ij} = f_i $ for $ i \leq j $ such 
that $ (f_i \circ \iota_i^*) $ converges strictly, then there exists a unique multiplier morphism $ f: \varinjlim_{i \in I} X_i \rightarrow Y $ such 
that $ f \circ \iota_i = f_i $ for all $ i \in I $. 
Using this mapping property one checks that a direct limit $ \varinjlim_{i \in I} X_i $ is unique up to isomorphism.  
However, not every inductive system in $ \BV $ need to admit a direct limit in $ \BV $. 

If $ \BV $ is a $ C^* $-category with finite direct sums and $ (X_i)_{i \in I} $ a family of objects in $ \BV $ then we obtain an inductive 
system $ (Y_F)_{F \in \F} $ over the directed set $ \F $ of finite subsets of $ I $ by considering $ Y_F = \bigoplus_{i \in \F} X_i \in \BV $ with the 
canonical inclusions as connecting maps. 
A direct limit of the system $ (Y_F)_{F \in \F} $ is nothing but a direct sum of the family $ (X_i)_{i \in I} $. 

The ind-category $ \Ind \BV $ of a $ C^* $-category $ \BV $ has the following closure property under inductive limits in the sense of 
Definition \ref{definductivelimit}. 

\begin{prop} 
Let $ \BV $ be a $ C^* $-category. Then every inductive system $ (X_i)_{i \in I} $ of objects in $ \BV $ admits an inductive limit in $ \Ind \BV $. 
Every countable inductive system of objects in $ \BV $ admits an inductive limit in $ \ind \BV $. 
\end{prop} 

\proof An inductive system $ (X_i)_{i \in I} $ in $ \BV $ can be viewed as an object $ X = (X_i)_{i \in I} \in \Ind \BV $, and we claim 
that $ X $ together with the canonical multiplier morphisms $ \iota_j: X_j \rightarrow X $ is an inductive limit 
of the inductive system $ (X_i)_{i \in I} $ of constant ind-objects in $ \Ind \BV $. 
In fact, we have $ \iota_j \circ \iota_{ij} = \iota_i $ for $ i \leq j $ 
and $ \iota_j \circ \iota_j^* \rightarrow \id $ in the strict topology by the definition of the morphism spaces in $ \Ind \BV $. 
The assertion for the countable ind-category $ \ind \BV $ is obtained in the same way. \qed 

We will be mainly interested in countably additive $ C^* $-categories. If $ \BV \in C^* \Lin $ then every countable inductive system in $ \BV $, 
viewed as an object of $ \Ind \BV $, is isomorphic to an inductive system of the form $ (V_n)_{n \in \mathbb{N}} $ 
with $ V_n = \bigoplus_{j = 1}^n X_j $ for some family $ (X_j)_{j \in \mathbb{N}} $ of objects of $ \BV $. 
It follows that every category in $ C^* \Lin $ is automatically closed under countable inductive limits. In particular, the canonical nondegenerate 
linear $ * $-functor $ \BV \rightarrow \ind \BV $ into the countable ind-category is an equivalence for $ \BV \in C^* \Lin $.

\subsection{Finitely accessible $ C^* $-categories} 

The theory of accessible and locally presentable categories is concerned with categories built using inductive limits from a set of presentable 
objects \cite{ARlocallypresentable}. This does not apply directly to $ C^* $-categories, 
but an analogous setup can be devised as follows. 

Let $ \BV \in C^* \Lin $ be a countably additive $ C^* $-category. We shall say that an object $ P \in \BV $ is \emph{finitely presentable} if $ \BV(P,P) $ is unital. 
In this case we have 
$$
\varinjlim_i \BV(P, V_i) \cong \BV(P, \varinjlim_i V_i) 
$$
for all countable inductive systems $ (V_i)_{i \in I} $.  
We let $ P(\BV) $ be the full subcategory of $ \BV $ consisting of all finitely presentable objects. The $ C^* $-category $ P(\BV) $ is closed under 
finite direct sums and subobjects. Note that any nondegenerate linear $ * $-functor $ \Bf: \BV \rightarrow \BW $ restricts to a unital 
linear $ * $-functor $ P(\BV) \rightarrow P(\BW) $. 

\begin{definition} \label{defcstaraccessible}
We say that a countably additive $ C^* $-category $ \BV $ is finitely accessible if there exists a set $ \P \subset \BV $ of finitely presentable objects such that 
every object in $ \BV $ is isomorphic to a subobject of a direct limit of some countable inductive system with objects from $ \P $. 
\end{definition} 

Finitely accessible $ C^* $-categories can equivalently be described as ind-categories of finitely linear $ C^* $-categories as follows. 

\begin{prop} \label{accessiblechar}
Let $ \BV \in C^* \Lin $ be a finitely accessible $ C^* $-category and let $ P(\BV) $ be the full subcategory of all finitely presentable objects in $ \BV $. 
Then $ \BV \simeq \ind P(\BV) $. 
\end{prop} 

\proof The canonical inclusion functor $ P(\BV) \rightarrow \BV $ 
induces a nondegenerate linear $ * $-functor $ \Bi: \ind P(\BV) \rightarrow \BV $ since $ \BV $ is closed under countable direct sums. 
From the definition of finite accessibility if follows that $ \Bi $ is essentially surjective. Moreover, $ \Bi $ is fully faithful on $ P(\BV) $ by 
construction, and it is also fully faithful on all countable direct sums of objects from $ P(\BV) $ in $ \ind P(\BV) $. 
Since every object of $ \ind P(\BV) $ is isomorphic to such a direct sum the claim follows. \qed 

Note in particular that the category $ \Hilb_A $ of countably generated Hilbert modules over 
a unital $ C^* $-algebra $ A $ is finitely accessible, and that $ P(\Hilb_A) $ identifies with the category $ \Hilb_A^f $ of finitely generated projective 
Hilbert $ A $-modules. 

In general, finitely presentable objects in countably additive $ C^* $-categories may be scarce. Even a basic example like the category $ \Hilb_A $ 
for $ A = C_0(\mathbb{R}) $ contains no nonzero finitely presentable objects. 
For such categories it is not obvious how to give intrinsic meaning to accessibility, and it seems best not to try to approximate them from smaller subcategories. 

Note also that while every object in the multiplier category $ M\BV $ of a $ C^* $-category $ \BV \in C^* \Lin $ is finitely presentable with respect to $ M\BV $, 
multiplier categories of nonzero $ C^* $-categories are not closed under countable direct sums, so that such categories are not finitely accessible 
in the sense of Definition \ref{defcstaraccessible}.

\subsection{Direct products of $ C^* $-categories} 

Let $ I $ be an index set and let $ (\BV_i)_{i \in I} $ be a family of $ C^* $-categories. 
The \emph{direct product} $ \prod_{i \in I} \BV_i $ is the category whose objects are families $ (X_i)_{i \in I} $ of objects $ X_i \in \BV_i $, with morphisms
$$
\bigg(\prod_{i \in I} \BV_i\bigg)((X_i), (Y_i)) = \bigoplus_{i \in I} \BV_i(X_i, Y_i). 
$$
Here the direct sum on the right hand side is taken in the $ c_0 $-sense. Equipped with the supremum norm and entrywise operations on morphisms this 
becomes a $ C^* $-category. 
The appearance of a direct sum instead of a direct product in the above formula may be surprising at first sight, but this choice is crucial for the 
validity of some of the arguments below. 

We obtain canonical full and essentially surjective nondegenerate linear $ * $-functors $ \pi_j: \prod_{i \in I} \BV_i \rightarrow \BV_j $ for all $ j \in I $. 
On the level of multipliers these functors induce canonical isomorphisms
$$
M\bigg(\prod_{i \in I} \BV_i\bigg)((X_i), (Y_i)) \cong \prod_{i \in I} M\BV_i(X_i, Y_i)
$$
for all objects $ (X_i), (Y_i) \in \prod_{i \in I} \BV_i $, with the product on the right hand side taken in the $ l^\infty $-sense. 

It follows in particular that $ \prod_{i \in I} \BV_i $ is subobject complete iff all $ \BV_i $ are subobject complete. 
If all categories $ \BV_i $ are contained in $ C^* \Lin $ then $ \prod_{i \in I} \BV_i $ is again contained in $ C^* \Lin $, with direct sums formed entrywise. 

\begin{prop} \label{categoricalcstardirectproducts}
Let $ I $ be an index set and let $ (\BV_i)_{i \in I} $ be a family of $ C^* $-categories in $ C^* \Lin $. 
For every $ \BU \in C^* \Lin $, postcomposition with the linear $ * $-functors $ \pi_j $ induces an equivalence of $ C^* $-categories 
$$
C^* \Lin(\BU, \prod_{i \in I} \BV_i) \simeq M \prod_{i \in I} C^*\Lin(\BU, \BV_i), 
$$
pseudonatural in $ \BU $. 
\end{prop} 

\proof We point out that the product $ \prod_{i \in I} C^*\Lin(\BU, \BV_i) $ on the right hand side is to be understood 
in the sense of $ C^* $-categories as defined above. 

Let us write $ F: C^* \Lin(\BU, \prod_{i \in I} \BV_i) \rightarrow M \prod_{i \in I} C^*\Lin(\BU, \BV_i) $ for the 
linear $ * $-functor given by $ F(\Bf) = (\pi_i \circ \Bf)_{i \in I} $ on objects and $ F(\sigma) = (\pi_i * \sigma)_{i \in I} $ on morphisms. 
If $ \Bh_i: \BU \rightarrow M\BV_i $ are nondegenerate linear $ * $-functors for $ i \in I $ then $ \Bh: \BU \rightarrow M \prod_{i \in I} \BV_i $ defined 
by $ \Bh(U) = (\Bh_i(U))_{i \in I} $ on objects and $ \Bh(f) = (\Bh_i(f))_{i \in I} $ on morphisms is a nondegenerate 
linear $ * $-functor with $ F(\Bh) = (\Bh_i)_{i \in I} $. Hence $ F $ is essentially surjective. 
If $ (\sigma_i)_{i \in I} $ is a bounded family of multiplier natural transformations $ \sigma_i: \pi_i \circ \Bh \Rightarrow \pi_i \circ \Bk $ for nondegenerate 
linear $ * $-functors $ \Bh, \Bk: \BU \rightarrow M\prod_{i \in I} \BV_i $, then assembling them pointwise yields a multiplier natural 
transformation $ \sigma: \Bh \Rightarrow \Bk $ such that $ F(\sigma) = (\sigma_i)_{i \in I} $. This means that $ F $ is full. 
Finally, note that $ F(\sigma) = 0 $ for $ \sigma: \Bh \Rightarrow \Bk $ implies $ (\pi_i * \sigma)(U) = 0 $ for $ U \in \BU $ and all $ i \in I $, 
which means that $ \sigma $ is zero. Hence $ F $ is faithful. \qed  

We note that the assertion of Proposition \ref{categoricalcstardirectproducts} fails if one does not allow multipliers in the definition of $ 1 $-morphisms 
in $ C^* \Lin $.

\section{Tensor products} \label{sectensorproducts}

In this section we review the construction of minimal and maximal tensor products of $ C^* $-categories, compare \cite{Mitchenercstarcategories}, 
\cite{DellAmbrogiounitarymodel}, and explain how to extend this to additive $ C^* $-categories. 

\subsection{Tensor products of $ C^* $-categories} 

If $ \BV_1, \BV_2 $ are $ C^* $-categories then their \emph{algebraic tensor product} $ \BV_1 \boxdot \BV_2 $ is the $ * $-category with objects all 
pairs $ (V_1, V_2) $ for $ V_1 \in \BV_1, V_2 \in \BV_2 $ and morphism spaces 
$$
(\BV_1 \boxdot \BV_2)((V_1, V_2),(W_1, W_2)) = \BV_1(V_1, W_1) \odot \BV_2(V_2, W_2), 
$$
where $ \odot $ denotes the algebraic tensor product. 
We will also write $ V_1 \boxdot V_2 = (V_1, V_2) $ to denote the objects in $ \BV_1 \boxdot \BV_2 $. 

The \emph{minimal tensor product} $ \BV_1 \boxdot_{\min} \BV_2 $ of $ \BV_1, \BV_2 $ has the same objects as $ \BV_1 \boxdot \BV_2 $, and the morphism spaces 
obtained by considering realisations $ \iota_1: \BV_1 \rightarrow M(\Hilb), \iota_2: \BV_2 \rightarrow M(\HILB) $, and the associated embeddings 
\begin{align*}
(\BV_1 \boxdot \BV_2)((V_1, V_2),(W_1, W_2)) &= \BV_1(V_1, W_1) \odot \BV_2(V_2, W_2) \\
&\subset \LH(\iota_1(V_1), \iota_1(W_1)) \otimes_{\min} \LH(\iota_1(V_2), \iota_2(W_2)) \\
&\subset \LH(\iota_1(V_1) \otimes \iota_2(V_2), \iota_1(W_1) \otimes \iota_2(W_2)) 
\end{align*}
for $ (V_1,V_2), (W_1, W_2) \in \BV_1 \boxdot \BV_2 $. More precisely, the morphism space in the minimal tensor product is
\begin{align*} 
(\BV_1 \boxdot_{\min} \BV_2)((V_1, V_2),(W_1, W_2)) &= [(\BV_1 \boxdot \BV_2)((V_1, V_2),(W_1, W_2))] \\
&\subset \LH(\iota_1(V_1) \otimes \iota_2(V_2), \iota_1(W_1) \otimes \iota_2(W_2)),
\end{align*} 
that is, the closure of the morphism space in the algebraic tensor product via this embedding. This is independent of the 
choice of $ \iota_1 $ and $ \iota_2 $. 
Using associativity of the minimal tensor product of $ C^* $-algebras one checks that the minimal tensor product $ \boxdot_{\min} $ is associative 
in a natural way. We will also write $ V_1 \boxdot_{\min} V_2 = (V_1, V_2) $ for objects in $ \BV_1 \boxdot_{\min} \BV_2 $. 

The \emph{maximal tensor product} $ \BV_1 \boxdot_{\max} \BV_2 $ of $ \BV_1, \BV_2 $ has again the same objects as $ \BV_1 \boxdot \BV_2 $, and morphism spaces 
obtained by taking the maximal completion of the morphism spaces in $ \BV_1 \boxdot \BV_2 $ as follows. 
Firstly, note that for any finite family of objects $ X_1, \dots, X_n $ in a $ C^* $-category $ \BX $ we obtain a canonical $ C^* $-algebra structure on 
$$
\bigoplus_{i,j = 1}^n \BX(X_i, X_j) = \BX^{\oplus}\bigg(\bigoplus_{i = 1}^n X_i, \bigoplus_{i = 1}^n X_i\bigg),  
$$
by considering the finite direct sum $ X_1 \oplus \cdots \oplus X_n $ in the finite additive completion $ \BX^{\oplus} $ of $ \BX $. 
Moreover, for $ m \leq n $ the obvious inclusion $ \bigoplus_{i,j = 1}^m \BW(X_i, X_j) \rightarrow \bigoplus_{i,j = 1}^n \BW(X_i, X_j) $ 
yields a hereditary subalgebra. 
Now given a pair of objects $ (V_1, V_2), (W_1, W_2) \in \BV_1 \boxdot \BV_2 $ let us choose finite families of 
objects $ (V_i)_{i \in I} \in \BV_1^{\oplus}, (W_j)_{j \in J} \in \BV_2^{\oplus} $ containing $ V_1, V_2 $ and $ W_1, W_2 $, respectively. 
We may then define the morphism space $ (\BV_1 \boxdot_{\max} \BV_2)((V_1, V_2),(W_1, W_2)) $ as the direct summand corresponding to $ (V_1, V_2), (W_1, W_2) $ in 
$$
\BV_1^{\oplus}\bigg(\bigoplus_{i \in I} V_i, \bigoplus_{i \in I} V_i\bigg) \otimes_{\max} \BV_2^{\oplus}\bigg(\bigoplus_{j \in J} W_j , \bigoplus_{j \in J} W_j \bigg). 
$$
Since the maximal tensor product of $ C^* $-algebras is compatible with inclusions of hereditary subalgebras, compare chapter 3 in \cite{BObook}, this 
does not depend on the choice of the finite families $ (V_i)_{i \in I}, (W_j)_{j \in J} $, and yields a unique $ C^* $-category structure. 
In particular, for $ V_1 = W_1, V_2 = W_2 $ we obtain 
$$
(\BV_1 \boxdot_{\max} \BV_2)((V_1, V_2), (V_1, V_2)) = \BV_1(V_1, V_1) \otimes_{\max} \BV_2(V_2, V_2), 
$$
and we have a canonical inclusion 
$$ 
(\BV_1 \boxdot_{\max} \BV_2)((V_1, V_2),(W_1, W_2)) \subset \BV_1(V_1 \oplus W_1, V_1 \oplus W_1) \otimes_{\max} \BV_2(V_2 \oplus W_2, V_2 \oplus W_2)
$$
for all objects $ (V_1, V_2), (W_1, W_2) $. 

The morphism space $ (\BV_1 \boxdot_{\max} \BV_2)((V_1, V_2),(W_1, W_2)) $ can be equivalently described 
as the completion of $ (\BV_1 \boxdot \BV_2)((V_1, V_2),(W_1, W_2)) $ with respect to 
$$
\|f\|_{\max} = \sup_{\pi: \BV_1 \boxdot \BV_2 \rightarrow M\HILB} \|\pi(f) \|,  
$$
where $ \pi $ runs over all nondegenerate linear $ * $-functors from $ \BV_1 \boxdot \BV_2 $ into $ M\HILB $. 
For unital $ C^* $-categories this is precisely the definition given in \cite{DellAmbrogiounitarymodel}. 
In the same way as before we will also write $ V_1 \boxdot_{\max} V_2 = (V_1, V_2) $ for objects in $ \BV_1 \boxdot_{\max} \BV_2 $. 
Using associativity of the maximal tensor product of $ C^* $-algebras one checks that the maximal tensor product $ \boxdot_{\max} $ is associative in 
a natural way. 

Let $ \BV_1, \BV_2, \BW $ be $ C^* $-categories. By a \emph{bilinear $ * $-functor} from $ \BV_1 \times \BV_2 $ to $ M\BW $ we mean a $ * $-functor 
$ \Bb: \BV_1 \times \BV_2 \rightarrow M\BW $ such that the maps 
$$ 
\Bb: \BV_1(V_1, W_1) \times \BV_2(V_2, W_2) \rightarrow M\BW(\Bb(V_1, W_1), \Bb(V_2, W_2)) 
$$ 
are bilinear for all $ V_i, W_i \in \BV_i $. 
Such a functor $ \Bb $ corresponds uniquely to a linear $ * $-functor $ \BV_1 \boxdot \BV_2 \rightarrow M\BW $, determined by the second arrow 
in the canonical factorisation 
$$
\BV_1(V_1, W_1) \times \BV_2(V_2, W_2) \rightarrow \BV_1(V_1, W_1) \odot \BV_2(V_2, W_2) \rightarrow M\BW(\Bb(V_1, W_1), \Bb(V_2, W_2)) 
$$
induced by $ \Bb $. 
This functor will be denoted $ L \Bb $, and referred to as the linearisation of $ \Bb $. We shall say that a 
bilinear $ * $-functor $ \Bb: \BV_1 \times \BV_2 \rightarrow M\BW $ is nondegenerate iff its linearisation 
$ L \Bb: \BV_1 \boxdot \BV_2 \rightarrow M\BW $ is.  

Observe that the definition of minimal and maximal tensor products yields canonical nondegenerate 
bilinear $ * $-functors $ \boxdot_{\min}: \BV_1 \times \BV_2 \rightarrow M(\BV_1 \boxdot_{\min} \BV_2) $ 
and $ \boxdot_{\max}: \BV_1 \times \BV_2 \rightarrow M(\BV_1 \boxdot_{\max} \BV_2) $. 
By the construction of the maximal tensor product, any nondegenerate bilinear $ * $-functor $ \Bb: \BV_1 \times \BV_2 \rightarrow M\BW $ 
prolongs uniquely to a nondegenerate linear $ * $-functor $ \BV_1 \boxdot_{\max} \BV_2 \rightarrow M\BW $. By slight abuse of 
notation, this functor will again be denoted by $ L \Bb $, and also referred to as the linearisation of $ \Bb $. 

\begin{lemma} \label{componentfunctors}
Let $ \BV_1, \BV_2, \BW \in C^* \Lin $ be countably additive $ C^* $-categories. 
Moreover assume that $ \Bb: \BV_1 \times \BV_2 \rightarrow M\BW $ is 
a nondegenerate bilinear $ * $-functor and let $ V_1 \in \BV_1, V_2 \in \BV_2 $. 
\begin{bnum} 
\item[a)] There exists a uniquely determined nondegenerate linear $ * $-functor $ \Bb(V_1, -): \BV_2 \rightarrow M\BW $ 
such that $ \Bb(V_1,-)(Y) = \Bb(V_1, Y) $ for all $ Y \in \BV_2 $, and 
$$
\Bb(f, g \circ h) = \Bb(V_1, g) \circ \Bb(f, h) = \Bb(f, g) \circ \Bb(V_1, h) 
$$
for all $ f \in \BV_1(V_1, V_1), g, h \in \BV_2(Y, Y) $. 
\item[b)] There exists a uniquely determined nondegenerate linear $ * $-functor $ \Bb(-, V_2): \BV_1 \rightarrow M\BW $ 
such that $ \Bb(-, V_2)(X) = \Bb(X, V_2) $ for all $ X \in \BV_1 $, and 
$$
\Bb(g \circ h, f) = \Bb(g, V_2) \circ \Bb(h, f) = \Bb(g, f) \circ \Bb(h, V_2) 
$$
for all $ f \in \BV_1(V_2, V_2), g, h \in \BV_1(X, X) $. 
\end{bnum} 
\end{lemma} 

\proof We shall only consider $ a) $, the case of $ b) $ being analogous. Note that we write $ \Bb(V_1, h) = \Bb(V_1,-)(h) $ for morphisms $ h $ in $ \BV_2 $.  

For each $ Y \in \BV_2 $ the linearisation of $ \Bb $ induces a 
unital $ * $-homomorphism $ L\Bb: M(\BV_1(V_1, V_1) \otimes_{\max} \BV_2(Y,Y)) \rightarrow M\BW(\Bb(V_1, Y), \Bb(V_1, Y)) $, and 
we define $ \Bb(V_1, -): \BV_2(Y,Y) \rightarrow M\BW(\Bb(V_1, Y), \Bb(V_1, Y)) $ by $ \Bb(V_1, g) = L\Bb(\id_{V_1} \odot g) $. 
We immediately get 
$$
\Bb(f, g \circ h) = \Bb(V_1, g) \circ \Bb(f, h) = \Bb(f, g) \circ \Bb(V_1, h) 
$$
for all $ f \in \BV_1(V_1, V_1), g, h \in \BV_2(Y, Y) $. 
To construct the action of $ \Bb(V_1,-) $ on morphisms $ g \in \BV_2(X,Z) $ consider the same construction for $ Y = X \oplus Z $ and 
identify $ \BV_2(X,Z) $ as a corner in $ \BV_2(Y,Y) $. One checks that this yields a well-defined nondegenerate linear $ * $-functor $ \Bb(V_1,-) $ as claimed. 
Uniqueness of $ \Bb(V_1, -) $ on objects is clear, and uniqueness on morphisms follows from nondegeneracy of $ \Bb $ and the defining formulas. \qed 

Lemma \ref{componentfunctors} implies in particular that a nondegenerate bilinear $ * $-functor preserves direct sums in both variables in the obvious sense. 

Let us state the basic functoriality properties of minimal and maximal tensor products. Assume that $ \BV_1, \BV_2, \BW_1, \BW_2 $ are $ C^* $-categories and 
let $ \Bf: \BV_1 \rightarrow \BW_1, \Bg: \BV_2 \rightarrow \BW_2 $ be linear $ * $-functors. 
From the definition of the algebraic tensor product $ \boxdot $ it is obvious that we obtain an induced 
linear $ * $-functor $ \Bf \boxdot \Bg: \BV_1 \boxdot \BV_2 \rightarrow \BW_1 \boxdot \BW_2 $ such that $ (\Bf \boxdot \Bg)(V_1, V_2) = (\Bf(V_1), \Bg(V_2)) $ 
on objects, and acting by $ (\Bf \boxdot \Bg)(f \odot g) = \Bf(f) \odot \Bg(g) $ on morphism spaces. 
This extends canonically to linear $ * $-functors $ \Bf \boxdot_{\max} \Bg: \BV_1 \boxdot_{\max} \BW_2 \rightarrow \BW_1 \boxdot_{\max} \BW_2 $ 
and $ \Bf \boxdot_{\min} \Bg: \BV_1 \boxdot_{\min} \BW_2 \rightarrow \BW_1 \boxdot_{\min} \BW_2 $, respectively.

\subsection{Tensor products of additive $ C^* $-categories} 

For finitely additive or countably additive $ C^* $-categories the above constructions have to be adapted in order to admit direct sums and subobjects. 

Let $ \BV_1, \BV_2 $ be (finitely/countably) additive $ C^* $-categories. We define $ \BV_1 \boxtimes_{\min} \BV_2 $ to be the subobject completion 
of $ \BV_1 \boxdot_{\min} \BV_2 $, and refer to it again as the minimal tensor product of $ \BV_1 $ and $ \BV_2 $. Similarly, the maximal tensor 
product $ \BV_1 \boxtimes_{\max} \BV_2 $ is the subobject completion of $ \BV_1 \boxdot_{\max} \BV_2 $. 
For objects $ V_1 \in \BV_1, V_2 \in \BV_2 $ we will write $ V_1 \boxtimes_{\min} V_2 $ and $ V_1 \boxtimes_{\max} V_2 $, respectively, for 
$ V_1 \boxdot_{\min} V_2 $ and $ V_1 \boxdot_{\max} V_2 $ viewed as objects in the subobject completions.  

The following argument shows that $ \boxtimes_{\min} $ and $ \boxtimes_{\max} $ both preserve the class of finitely (countably) 
additive $ C^* $-categories, respectively. 

\begin{prop} 
Let $ \BV_1, \BV_2 $ be finitely (countably) additive $ C^* $-categories. Then the minimal and maximal tensor 
products $ \BV_1 \boxtimes_{\min} \BV_2, \BV_1 \boxtimes_{\max} \BV_2 $ are again finitely (countably) additive. 
\end{prop} 

\proof Consider the case of the maximal tensor product $ \BV_1 \boxtimes_{\max} \BV_2 $ for countably additive categories $ \BV_1, \BV_2 \in C^* \Lin $. 
Since $ \BV_1 \boxtimes_{\max} \BV_2 $ is subobject complete by construction it suffices to check that this category admits countable 
direct sums. 

If $ (X_n)_{n \in \mathbb{N}} $ is a family of objects in $ \BV_1 $ and $ (Y_n)_{n \in \mathbb{N}} $ a family of objects in $ \BV_2 $ 
then $ (\bigoplus_{m \in \mathbb{N}} X_m) \boxtimes_{\max} (\bigoplus_{n \in \mathbb{N}} Y_n) $ is a direct sum of the 
family $ (X_m \boxtimes_{\max} Y_n)_{m,n \in \mathbb{N}} $ in $ \BV_1 \boxtimes_{\max} \BV_2 $. 
Now if $ (Z_n)_{n \in \mathbb{N}} $ is an arbitrary countable family of objects in $ \BV_1 \boxtimes_{\max} \BV_2 $ 
then for every $ n \in \mathbb{N} $ we find $ X_n \in \BV_1, Y_n \in \BV_2 $ such that $ Z_n $ is a subobject of $ X_n \boxtimes_{\max} Y_n $. 
Then we can realise the direct sum $ \bigoplus_{n \in \mathbb{N}} Z_n $ as subobject 
of $ (\bigoplus_{n \in \mathbb{N}} X_n) \boxtimes_{\max} (\bigoplus_{n \in \mathbb{N}} Y_n) $. 

The proof for minimal tensor products is analogous, and in the case of finitely additive categories it suffices to consider only finite direct sums instead. \qed 

Let $ \BV_1, \BV_2, \BW \in C^* \Lin $ and let $ \Bb, \Bc: \BV_1 \times \BV_2 \rightarrow M\BW $ be nondegenerate bilinear $ * $-functors. 
A multiplier natural transformation $ \phi: \Bb \Rightarrow \Bc $ is a natural transformation of the underlying $ * $-functors 
such that $ \phi(V_1, V_2) \in M(\Bb(V_1, V_2), \Bc(V_1, V_2)) $ is uniformly bounded for all $ (V_1, V_2) \in \BV_1 \times \BV_2 $. 
The collection of all nondegenerate bilinear $ * $-functors from $ \BV_1 \times \BV_2 $ to $ M\BW $ together with their multiplier natural transformations 
forms a $ C^* $-category which we will denote by $ C^* \Bilin(\BV_1, \BV_2; \BW) $. 
 
By construction, we have canonical nondegenerate bilinear $ * $-functors $ \boxtimes_{\min}: \BV_1 \times \BV_2 \rightarrow \BV_1 \boxtimes_{\min} \BV_2 $ 
and $ \boxtimes_{\max}: \BV_1 \times \BV_2 \rightarrow \BV_1 \boxtimes_{\max} \BV_2 $. On the level of objects these functors map $ (V_1,V_2) $ 
to $ V_1 \boxtimes_{\min} V_2 $ and $ V_1 \boxtimes_{\max} V_2 $, respectively. 

\begin{prop} \label{universalpropmaxtensor}
Let $ \BV_1, \BV_2, \BW \in C^* \Lin $. Then precomposition with the canonical 
bilinear $ * $-functor $ \boxtimes_{\max}: \BV_1 \times \BV_2 \rightarrow \BV_1 \boxtimes_{\max} \BV_2 $ induces an equivalence 
$$
C^* \Bilin(\BV_1, \BV_2; \BW) \simeq C^* \Lin(\BV_1 \boxtimes_{\max} \BV_2, \BW) 
$$
of $ C^* $-categories, pseudonatural in $ \BW $. 
\end{prop} 

\proof Consider the functor $ F: C^* \Lin(\BV_1 \boxtimes_{\max} \BV_2, \BW) \rightarrow C^* \Bilin(\BV_1, \BV_2; \BW) $ given 
by $ F(\Bf) = \Bf \circ \boxtimes_{\max}, F(\sigma) = \sigma * \boxtimes_{\max} $. 
It is easy to check that $ F $ is a linear $ * $-functor, and by definition of $ \BV_1 \boxtimes_{\max} \BV_2 $ this functor is essentially surjective. 

Assume that $ \sigma: \Bf \Rightarrow \Bg $ in $ C^* \Lin(\BV_1 \boxtimes_{\max} \BV_2, \BW) $ satisfies $ F(\sigma) = 0 $. 
Then $ \sigma(V_1 \boxtimes_{\max} V_2): \Bf(V_1 \boxtimes_{\max} V_2) \rightarrow \Bg(V_1 \boxtimes_{\max} V_2) $ vanishes 
for all $ V_1 \in \BV_1, V_2 \in \BV_2 $, and hence also on all subobjects of such tensor products by naturality. 
We conclude $ \sigma = 0 $, which means that $ F $ is faithful. 
Conversely, let $ \sigma: \Bf \circ \boxtimes_{\max} \Rightarrow \Bg \circ \boxtimes_{\max} $ in $ C^* \Bilin(\BV_1, \BV_2; \BW) $ be given. 
Then $ \sigma $ can be viewed as multiplier natural transformation between the corresponding linear $ * $-functors defined on $ \BV_1 \boxdot_{\max} \BV_2 $, 
and extends canonically from $ \BV_1 \boxdot_{\max} \BV_2 $ to its subobject completion. It follows that $ F $ is full. \qed 

Let us next record the functoriality properties of tensor products in the additive setting. Let $ \BV_1, \BV_2, \BW_1, \BW_2 \in C^* \Lin $ and 
let $ \Bf: \BV_1 \rightarrow M\BW_1, \Bg: \BV_2 \rightarrow M\BW_2 $ be nondegenerate linear $ * $-functors. 
Then we obtain induced nondegenerate linear $ * $-functors $ \Bf \boxtimes_{\max} \Bg: \BV_1 \boxtimes_{\max} \BW_2 \rightarrow M(\BW_1 \boxtimes_{\max} \BW_2) $ 
and $ \Bf \boxtimes_{\min} \Bg: \BV_1 \boxtimes_{\min} \BW_2 \rightarrow M(\BW_1 \boxtimes_{\min} \BW_2) $ which 
map $ (V_1, V_2) $ to $ (\Bf(V_1), \Bg(V_2)) $ and act by sending $ f \boxdot g $ to $ \Bf(f) \odot \Bg(g) $ on morphism spaces.

\begin{theorem} \label{symmon2cat}
Both the maximal and the minimal tensor product determine symmetric monoidal structures on the $ 2 $-category $ C^* \Lin $. 
\end{theorem} 

\proof For the axioms of symmetric monoidal bicategories see \cite{GPStricategories}, \cite{McCruddenbalanced}, \cite{Schommerpriesthesis}. 
Since the arguments for minimal and maximal tensor products are analogous we shall restrict attention to $ \boxtimes_{\max} $. 

If $ \BV_1, \BV_2, \BV_3 $ are in $ C^* \Lin $ then we have obvious associativity 
equivalences $ (\BV_1 \boxtimes_{\max} \BV_2) \boxtimes_{\max} \BV_3 \rightarrow \BV_1 \boxtimes_{\max} (\BV_2 \boxtimes_{\max} \BV_3) $ 
sending $ (V_1 \boxtimes_{\max} V_2) \boxtimes_{\max} V_3 $ to $ V_1 \boxtimes_{\max} (V_2 \boxtimes_{\max} V_3) $ for objects $ V_i \in \BV_i $. 
The unit object in $ C^* \Lin $ is $ 1 = \Hilb = \Hilb_\mathbb{C} $, and on the level of objects the left and right unitor 
equivalences $ l: \Hilb \boxtimes_{\max} \BV \rightarrow \BV, r: \BV \rightarrow \BV \boxtimes_{\max} \Hilb $ map $ \mathbb{C} \boxtimes_{\max} V $ to $ V $ 
and $ V $ to $ V \boxtimes_{\max} \mathbb{C} $, respectively.  
It is straightforward to check that this data satisfies the axioms of a monoidal bicategory. 

The symmetry equivalence $ \sigma: \BV_1 \boxtimes_{\max} \BV_2 \rightarrow \BV_2 \boxtimes_{\max} \BV_1 $ for $ \BV_1, \BV_2 \in C^* \Lin $ is 
given by $ \sigma(V_1 \boxtimes_{\max} V_2) = V_2 \boxtimes_{\max} V_1 $ on the level of objects. Despite the fact that the axioms for a symmetric monoidal 
bicategory are rather unwieldy, in the case at hand they boil down to elementary properties of the maximal tensor product of $ C^* $-algebras. \qed 

Maximal and minimal tensor products of countably additive $ C^* $-categories are compatible with exterior tensor products of Hilbert modules in 
the following sense. 

\begin{prop} \label{additivehilbert}
Let $ A, B $ be $ C^* $-algebras. Then there are canonical equivalences
$$ 
\Hilb_A \boxtimes_{\min} \Hilb_B \simeq \Hilb_{A \otimes_{\min} B}, \qquad \Hilb_A \boxtimes_{\max} \Hilb_B \simeq \Hilb_{A \otimes_{\max} B}
$$ 
of $ C^* $-categories. 
\end{prop} 

\proof The proof for minimal and maximal tensor products is similar, so let us only consider the case of maximal tensor products. 

The maximal exterior tensor product $ \E \otimes_{\max} \F $ of $ \E \in \Hilb_A, \F \in \Hilb_B $ is the completion of $ \E \odot \F $ with respect to 
the $ A \otimes_{\max} B $-valued inner product given by 
$$
\bra \xi \otimes \eta, \zeta \otimes \kappa \ket = \bra \xi, \zeta \ket \otimes \bra \eta, \kappa \ket.  
$$
This defines canonical fully faithful linear $ * $-functor $ \Hilb_A \boxtimes_{\max} \Hilb_B \rightarrow \Hilb_{A \otimes_{\max} B} $. 
Since every Hilbert module in $ \Hilb_{A \otimes_{\max} B} $ is isomorphic to a direct summand of the standard Hilbert module $ \HH_{A \otimes_{\max} B} $ 
it follows that this embedding is essentially surjective. \qed 

An analogue of Proposition \ref{additivehilbert} holds in the finitely additive setting as well, that is, we have equivalences 
$$ 
\Hilb^f_A \boxtimes_{\min} \Hilb^f_B \simeq \Hilb^f_{A \otimes_{\min} B}, \qquad \Hilb^f_A \boxtimes_{\max} \Hilb^f_B \simeq \Hilb^f_{A \otimes_{\max} B}
$$ 
for the tensor products of categories of finitely generated projective Hilbert modules over unital $ C^* $-algebras $ A,B $.

\subsection{Tensor products of finitely accessible $ C^* $-categories} 

In this section we discuss the special case of tensor products of finitely accessible $ C^* $-categories, and also relate $ C^* $-tensor products 
to the Deligne tensor product in the purely algebraic setting. 

\begin{prop} \label{finitetensorproducts}
If $ \BV_1, \BV_2 $ are finitely accessible $ C^* $-categories then $ \BV_1 \boxtimes_{\max} \BV_2 $ is finitely accessible and 
$$ 
P(\BV_1 \boxtimes_{\max} \BV_2) \simeq P(\BV_1) \boxtimes_{\max} P(\BV_2). 
$$
Similarly, if $ \BV_1, \BV_2 $ are finitely additive $ C^* $-categories then 
$$ 
\ind(\BV_1 \boxtimes_{\max} \BV_2) \simeq \ind(\BV_1) \boxtimes_{\max} \ind(\BV_2). 
$$
Both statements hold for minimal tensor products as well. 
\end{prop} 

\proof The arguments are analogous for maximal and minimal tensor products, so we shall only consider maximal tensor products 
and abbreviate $ \boxtimes = \boxtimes_{\max} $. 

Assume first that $ \BV_1, \BV_2 $ are finitely accessible $ C^* $-categories. The tensor product object $ V_1 \boxtimes V_2 $ of finitely presented 
objects $ V_i \in P(\BV_i) $ is clearly finitely presented in $ \BV_1 \boxtimes \BV_2 $. 
Since every object of $ \BV_1 \boxtimes \BV_2 $ can be written as a subobject of some countable 
inductive limit of such objects it follows that the category $ \BV_1 \boxtimes \BV_2 $ is finitely accessible. 

The above argument shows that we have a canonical fully faithful embedding $ P(\BV_1) \boxtimes P(\BV_2) \rightarrow P(\BV_1 \boxtimes \BV_2) $.  
In order to show that this embedding is essentially surjective let $ V \in \BV_1, W \in \BV_2 $ and 
write $ V \cong \varinjlim_{i \in \mathbb{N}} V_i, W \cong \varinjlim_{j \in \mathbb{N}} W_j $ as countable inductive limits of finitely presentable objects, 
and denote by $ \iota^V_i, \iota^W_j, \pi^V_i, \pi^W_j $ the associated isometries and their adjoints. 
If $ p \in (\BV_1 \boxtimes \BV_2)(V \boxtimes W, V \boxtimes W) $ is a projection and $ \epsilon > 0 $ then there 
exist $ m, n \in \mathbb{N} $ such that $ \|P_{mn} \circ p \circ P_{mn} - p \| < \epsilon $, 
where $ P_{mn} = (\iota^V_m \circ \pi^V_m) \boxtimes (\iota^W_n \circ \pi^W_n) $. 
Choosing $ \epsilon $ small enough and using functional calculus we obtain a projection $ q \in (\BV_1 \boxtimes \BV_2)(V_m \boxtimes W_n, V_m \boxtimes W_n) $ 
such that $ \|p - (\iota^V_m \boxtimes \iota^W_n) \circ q \circ (\pi^V_m \boxtimes \pi^W_n)\| < 1 $. 
The projections $ p $ and $ q $ are unitarily equivalent in $ M(\BV_1 \boxtimes \BV_2)(V \boxtimes W, V \boxtimes W) $, 
which means that the subobjects corresponding to them are isomorphic. It follows that every finitely presented object in $ \BV_1 \boxtimes \BV_2 $ 
arises as a subobject of the tensor product $ V_1 \boxtimes V_2 $ of finitely presented objects $ V_i \in P(\BV_i) $. 

Now assume that $ \BV_1, \BV_2 $ are finitely additive $ C^* $-categories. The obvious 
linear $ * $-functor $ \ind(\BV_1 \boxtimes \BV_2) \rightarrow \ind(\BV_1) \boxtimes \ind(\BV_2) $ 
is fully faithful and essentially surjective, hence an equivalence. \qed 

Assume that $ \BV_1, \BV_2 $ are finitely additive $ C^* $-categories such that all morphism spaces are finite dimensional. 
Then the categories $ \BV_1, \BV_2 $ are semisimple and in particular $ \mathbb{C} $-linear abelian.  
In this case there are no completions involved in the definition of the minimal or maximal tensor products, 
and the category $ \BV_1 \boxtimes_{\max} \BV_2 \cong \BV_1 \boxtimes_{\min} \BV_2 $ is again $ \mathbb{C} $-linear abelian. We record the following 
observation. 

\begin{prop} 
Let $ \BV_1, \BV_2 $ be finitely additive $ C^* $-categories with finite dimensional morphism spaces. Then we have a canonical equivalence 
$$
\BV_1 \boxtimes_{\max} \BV_2 \simeq \BV_1 \boxtimes \BV_2, 
$$
where $ \BV_1 \boxtimes \BV_2 $ denotes the Deligne tensor product of the $ \mathbb{C} $-linear abelian categories underlying $ \BV_1, \BV_2 $. 
\end{prop} 

\proof The canonical bilinear $ * $-functor $ \BV_1 \times \BV_2 \rightarrow \BV_1 \boxtimes_{\max} \BV_2 $ is right exact in each variable, 
and therefore induces a linear functor $ \gamma: \BV_1 \boxtimes \BV_2 \rightarrow \BV_1 \boxtimes_{\max} \BV_2 $. If $ (X_i)_{i \in I} $ are representatives 
of the isomorphism classes of simple objects in $ \BV_1 $ and $ (Y_j)_{j \in J} $ the same for $ \BV_2 $, then both $ \BV_1 \boxtimes \BV_2 $ 
and $ \BV_1 \boxtimes_{\max} \BV_2 $ are semisimple with isomorphism classes of simple objects given 
by $ (X_i \boxtimes Y_j)_{i \in I,j \in J} $ and $ (X_i \boxtimes_{\max} Y_j)_{i \in I, j \in J} $, respectively. It follows that $ \gamma $ is an 
equivalence of categories. \qed 

Given arbitrary $ \BV, \BW \in C^* \lin $ the $ C^* $-category $ C^* \lin(\BV, \BW) $ of unital linear $ * $-functors from $ \BV $ to $ \BW $ 
and their natural transformations has finite direct sums, obtained by taking pointwise direct sums in $ \BW $. 
The category $ C^* \lin(\BV, \BW) $ is also closed under subobjects, again taken pointwise. 
It follows that $ C^* \lin(\BV, \BW) $ is contained in $ C^* \lin $. 

The following result shows that this construction yields an internal Hom in $ C^* \lin $, compare \cite{DellAmbrogiounitarymodel}. 

\begin{prop} \label{lininternalhom}
Let $ \BU, \BV, \BW \in C^* \lin $ be finitely additive $ C^* $-categories. Then there exists an equivalence 
$$
C^* \lin(\BU \boxtimes_{\max} \BV, \BW) \simeq C^* \lin(\BU, C^* \lin(\BV, \BW)) 
$$
of $ C^* $-categories, pseudonatural in all variables. 
\end{prop} 

\proof Assume that $ \Bf: \BU \boxtimes_{\max} \BV \rightarrow \BW $ is a unital linear $ * $-functor. Then for every $ U \in \BU $ the 
restriction $ \Bf(U \boxtimes_{\max} -) $ defines a unital linear $ * $-functor $ \BV \rightarrow \BW $, 
and we define $ F: C^* \lin(\BU \boxtimes_{\max} \BV, \BW) \rightarrow C^* \lin(\BU, C^* \lin(\BV, \BW)) $ by $ F(\Bf)(U) = \Bf(U \boxtimes_{\max} -) $ on objects
and $ F(\Bf)(f) = \Bf(f \boxtimes_{\max} -) $ on morphisms. Then $ F(\Bf) $ is a linear $ * $-functor,  
and if $ \sigma: \Bf \Rightarrow \Bg $ is a natural transformation we obtain a natural transformation $ F(\sigma): F(\Bf) \Rightarrow F(\Bg) $ 
by setting $ F(\sigma)(U)(V) = \sigma(U \boxtimes_{\max} V) $. 
In this way $ F $ becomes a linear $ * $-functor, and it is routine to verify that it is an equivalence and pseudonatural in all variables. \qed 

We remark that Proposition \ref{lininternalhom} does not generalise to the $ 2 $-category $ C^* \Lin $ of countably additive $ C^* $-categories in any obvious way.

\section{Bicolimits of additive $ C^* $-categories} \label{secbicolim}

The main aim of this section is to show that the $ 2 $-category $ C^* \Lin $ of countably additive $ C^* $-categories is closed under conical bicolimits. 
We shall restrict ourselves to the unitary setting, which means that we will assume that all multiplier natural transformations appearing in 
homomorphisms and transformations are unitary. 
For simplicity we will also consider only bicolimits indexed by $ 1 $-categories, noting that the discussion below extends to arbitrary indexing $ 2 $-categories 
with minor modifications. 

Fix a small unital category $ I $. We define an $ I $-diagram $ \iota: I \rightarrow C^* \Lin $ to be given 
by $ C^* $-categories $ \iota(i) = \BV_i \in C^* \Lin $ for all objects $ i \in I $, nondegenerate 
linear $ * $-functors $ \iota(i \to j) = \iota_{i,j}(i \to j): \BV_i \rightarrow M\BV_j $ for all morphisms $ i \to j $ in $ I $, and unitary 
natural isomorphisms $ \iota_{i,j,k}: \iota(j \to k) \circ \iota(i \to j) \Rightarrow \iota(i \to k) $ for all pairs of composable 
morphisms $ i \to j, j \to k $ in $ I $, as well as unitary natural isomorphisms $ \iota_i: 1_{\iota(i)} \Rightarrow \iota(1_i) $, such that the diagrams 
\begin{center}
\begin{tikzcd}[column sep = 15mm]
			\iota(k \to l) \circ \iota(j \to k) \circ \iota(i \to j) 
			\arrow[d, swap, "\id_{\iota(k \to l)} * \iota_{i, j, k}"]
  		\arrow[rr,"\iota_{j, k, l} * \id_{\iota(i \to j)}"]
			&&
			\iota(j \to l) \circ \iota(i \to j)
			\arrow[d, "\iota_{i, j, l}"]
			\\
			\iota(k \to l) \circ \iota(i \to k) 
			\arrow[rr, swap, "\iota_{i, k, l}"]
      &&
			\iota(i \to l) 
\end{tikzcd}
\end{center}
and 
\begin{center}
\begin{tikzcd}
			\iota(i \to j) \circ 1_{\iota(i)}
    	\arrow[dd, swap, "\id_{\iota(i \to j)} * \iota_i"]
			\arrow[dr, "\id"]
			&&
			1_{\iota(j)} \circ \iota(i \to j) 
			\arrow[dl, swap, "\id"]
			\arrow[dd, "\iota_j * \id_{\iota(i \to j)}"]
			\\
			&
			\iota(i \to j)
			&
			\\
			\iota(i \to j) \circ \iota(1_i)
			\arrow[ur, swap, "\iota_{i,i,j}"]
			&&
			\iota(1_j) \circ \iota(i \to j)
      \arrow[ul, "\iota_{i,j,j}"]
\end{tikzcd}
\end{center}
commute. 
In the sequel we will sometimes use the notation $ (\BV_i)_{i \in I} $ for such an $ I $-diagram, suppressing the connecting functors and 
natural isomorphism data. 
Note that an $ I $-diagram defines a homomorphism $ \iota: I \rightarrow C^* \Lin $, where $ I $ is viewed as a $ 2 $-category with only identity $ 2 $-morphisms. 
This homomorphism has the extra property that all multiplier natural isomorphisms appearing in the definition are unitary. 

By definition, a transformation $ \Bf: \iota \rightarrow \eta $ of $ I $-diagrams $ \iota, \eta: I \rightarrow C^* \Lin $ consists of nondegenerate 
linear $ * $-functors $ \Bf_i: \iota(i) \rightarrow M\eta(i) $ for $ i \in I $ together with unitary 
natural isomorphisms $ \Bf_{i,j}(i \to j): \eta(i \to j) \circ \Bf_i \Rightarrow \Bf_j \circ \iota(i \to j) $ for all morphisms $ i \to j $ in $ I $ such that 
the diagrams 
\begin{center}
\begin{tikzcd}[column sep = -5mm]
		  &
			\eta(i \to j \to k) \circ \Bf_i
			\arrow[dr,"\Bf_{i, k}(i \to j \to k)"]
			&
			\\
			\eta(j \to k) \circ \eta(i \to j) \circ \Bf_i
			\arrow[ur, "\eta_{i, j, k} * \id_{\Bf_i}"]
			\arrow[d, swap, "\id_{\eta(j \to k)} * \Bf_{i,j}(i \to j)"]
			&&
			\Bf_k \circ \iota(i \to j \to k) 
			\\
			\eta(j \to k) \circ \Bf_j \circ \iota(i \to j) 
			\arrow[dr, swap, "\id"]
			&&
			\Bf_k \circ \iota(j \to k) \circ \iota(i \to j)   
			\arrow[u, swap, "\id_{\Bf_k} * \iota_{i,j,k}"]
			\\
			&
			\eta(j \to k) \circ \Bf_j \circ \iota(i \to j)
			\arrow[ur, swap, "\Bf_{j, k}(j \to k) * \id_{\iota(i \to j)}"]
			&
\end{tikzcd}
\end{center}
and 
\begin{center}
\begin{tikzcd}[column sep = -5mm]
			&&
			\qquad \Bf_i \qquad 
			\arrow[drr,"\id", start anchor = {[xshift = -3mm, yshift = 1mm]}, end anchor = {[xshift = 2mm]}]
			&&
			\\
			1_{\eta(i)} \circ \Bf_i
			\arrow[urr,"\id", start anchor = {[xshift = -2mm]}, end anchor = {[xshift = 3mm, yshift = 1mm]}]
			\arrow[dr, swap, "\eta_i * \id_{\Bf_i}"]
			&&&&
			\Bf_i \circ 1_{\iota(i)} 
			\arrow[dl, "\id_{\Bf_i} * \iota_i"]
			\\
			&
			\eta(1_i) \circ \Bf_i
			\arrow[rr, swap, "\Bf_{i,i}(1_i)"]
			&&
			\Bf_i \circ \iota(1_i) 
			&
\end{tikzcd}
\end{center}
commute. We will sometimes just write $ \Bf = (\Bf_i)_{i \in I}: (\BV_i)_{i \in I} \rightarrow (\BW_i)_{i \in I} $ for a transformation, suppressing the 
natural isomorphisms from the notation. 

A modification $ \Gamma: \Bf \Rightarrow \Bg $ between two transformations $ \Bf, \Bg: \iota \rightarrow \eta $ consists of a uniformly bounded family 
of multiplier natural transformations $ \Gamma_i: \Bf_i \Rightarrow \Bg_i $ for $ i \in I $ such that 
\begin{center}
\begin{tikzcd}[column sep = 15mm]
			\eta(i \to j) \circ \Bf_i
			\arrow[r, "\id_{\eta(i \to j)} * \Gamma_i"]
			\arrow[d, swap, "\Bf_{i, j}(i \to j)"]
			&
			\eta(i \to j) \circ \Bg_i
			\arrow[d, "\Bg_{i, j}(i \to j)"]
			\\
			\Bf_j \circ \iota(i \to j)
			\arrow[r, swap, "\Gamma_j * \id_{\iota(i \to j)}"]
			&
			\Bg_j \circ \iota(i \to j)
\end{tikzcd}
\end{center}
is commutative for every morphism $ i \to j $ in $ I $. It is called unitary if all $ \Gamma_i $ are unitary natural isomorphisms. 
We note that the category $ [I, C^* \Lin]((\BV_i)_{i \in I}, (\BW_i)_{i \in I}) $ of all 
transformations $ (\BV_i)_{i \in I} \rightarrow (\BW_i)_{i \in I} $ and their modifications is naturally a $ C^* $-category. 

For $ \BW \in C^* \Lin $ we obtain the constant $ I $-diagram $ \Delta(\BW) $ by setting $ \Delta(\BW)(i) = \BW $ and $ \Delta(\BW)(i \rightarrow j) = 1_\BW $ 
for all $ i \rightarrow j $ in $ I $. Writing $ [I, C^* \Lin] $ for the $ 2 $-category of all $ I $-diagrams in $ C^* \Lin $, this defines a strict homomorphism 
$ \Delta: C^* \Lin \rightarrow [I, C^* \Lin] $. 

\begin{definition}
Let $ (\BV_i)_{i \in I} $ be an $ I $-diagram and let $ \Bf = (\Bf_i)_{i \in I}: (\BV_i)_{i \in I} \rightarrow \Delta(\BX) $ be a transformation. We say 
that $ \Bf $ is cofibrant if $ \Bf_i(X) = \Bf_j(Y) $ for any $ i,j \in I $ and $ X \in \BV_i, Y \in \BV_j $ implies $ i = j $ and $ X = Y $. 
\end{definition} 

Given an arbitrary transformation $ \Bf = (\Bf_i)_{i \in I}: (\BV_i)_{i \in I} \rightarrow \Delta(\BX) $ let us construct a category $ CF(\BX) \in C^* \Lin $ 
and a cofibrant transformation $ CF(\Bf) = (CF(\Bf)_i)_{i \in I}: (\BV_i)_{i \in I} \rightarrow \Delta(CF(\BX)) $ as follows. 

Firstly, form the disjoint union $ \Lambda = \bigcup_{i \in I} \Ob(\BV_i) $ of all object sets of the categories $ \BV_i $, and let $ CF(\BX) $ be the 
countably additive $ C^* $-category with objects all pairs $ (\lambda, X) $ with $ \lambda \in \Lambda $ and $ X \in \Ob(\BX) $, morphism sets 
$$ 
CF(\BX)((\lambda, X), (\rho,Y)) = \BX(X,Y),  
$$ 
and the structures on morphisms such that the obvious forgetful map $ \phi_\BX: CF(\BX) \rightarrow \BX $, given by $ \phi_\BX(\lambda, X) = X $ 
on objects and acting as the identity on morphisms is a fully faithful linear $ * $-functor. 
By construction, the functor $ \phi_\BX $ is then an equivalence of $ C^* $-categories. 

Next define nondegenerate linear $ * $-functors $ CF(\Bf)_i: \BV_i \rightarrow MCF(\BX) $ by setting $  CF(\Bf)_i(V) = (V, \Bf_i(V)) $ for $ V \in \BV_i $ and 
$$ 
CF(\Bf)_i(f) = f \in MCF(\BX)((V, \Bf_i(V)), (W, \Bf_i(W))) = M\BX(\Bf_i(V), \Bf_i(W)) 
$$ 
for $ f \in \BV_i(V, W) $, and unitary natural isomorphisms $ CF(\Bf)_{i,j}(i \to j): CF(\Bf)_i \Rightarrow CF(\Bf)_j \circ \iota(i \to j) $ 
for every morphism $ i \to j $ in $ I $ by 
$$ 
CF(\Bf)_{i,j}(i \to j) = \Bf_{i,j}(i \to j). 
$$
This data yields a transformation $ CF(\Bf) $, since the extra variable just carries the source and target objects. 
We have $ \Delta(\phi_\BX) \circ CF(\Bf) = \Bf $ by construction, and note that $ \phi_\BX $ is pseudo-natural in $ \BX $. 

If $ \BW $ is a countably additive $ C^* $-category we write $ [I, C^* \Lin]_{cof}((\BV_i)_{i \in I}, \Delta(\BW)) $ for the full subcategory of 
$ [I, C^* \Lin]((\BV_i)_{i \in I}, \Delta(\BW)) $ consisting of all cofibrant transformations from $ (\BV_i)_{i \in I} $ to the constant diagram $ \Delta(\BW) $. 

\begin{lemma} \label{cofibrantlemma}
Let $ (\BV_i)_{i \in I} $ be an $ I $-diagram of countably additive $ C^* $-categories. Then postcomposition with $ \Delta(\phi_\BX) $ induces an equivalence 
of $ C^* $-categories 
$$ 
[I, C^* \Lin]_{cof}((\BV_i)_{i \in I}, \Delta(CF(\BX))) \rightarrow [I, C^* \Lin]((\BV_i)_{i \in I}, \Delta(\BX)), 
$$
pseudonatural in $ \BX $.  
\end{lemma} 

\proof For fixed $ \BX \in C^* \Lin $ we obtain a well-defined functor $ \Delta(\phi_\BX) \circ - $ from 
$ [I, C^* \Lin]_{cof}((\BV_i)_{i \in I}, \Delta(CF(\BX))) $ to $ [I, C^* \Lin]((\BV_i)_{i \in I}, \Delta(\BX)) $. 
Since $ \phi_\BX $ is an equivalence this functor is fully faithful, 
and our above arguments show that it is essentially surjective. Pseudonaturality is clear from the construction. \qed 

By the cardinality of a $ C^* $-category $ \BV \in C^* \Lin $ we shall mean the cardinality of the union all morphism spaces in $ \BV $. 
With these preparations in place we can now formulate and prove the following result. 

\begin{theorem} \label{bicolimitexistence}
Let $ I $ be a small category and let $ (\BV_i)_{i \in I} $ be an $ I $-diagram in $ C^* \Lin $. Then there exists a 
bicolimit $ \varinjlim_{i \in I} \BV_i \in C^* \Lin $, that is, there exists a transformation $ T: (\BV_i)_{i \in I} \rightarrow \Delta(\varinjlim_{i \in I} \BV_i) $
such that precomposition with $ T $ induces an equivalence of $ C^* $-categories
$$
C^* \Lin(\varinjlim_{i \in I} \BV_i, \BW) \simeq [I, C^* \Lin]((\BV_i)_{i \in I}, \Delta(\BW))
$$
for every $ \BW \in C^* \Lin $, pseudonatural in $ \BW $. 
\end{theorem} 

\proof Due to Lemma \ref{cofibrantlemma} it suffices to construct $ \varinjlim_{i \in I} \BV_i \in C^* \Lin $ and a 
transformation $ T: (\BV_i)_{i \in I} \rightarrow \Delta(\varinjlim_{i \in I} \BV_i) $ such that precomposition with $ T $ induces an equivalence 
of $ C^* $-categories 
$$
C^* \Lin(\varinjlim_{i \in I} \BV_i, \BW) \simeq [I, C^* \Lin]_{cof}((\BV_i)_{i \in I}, \Delta(\BW))
$$
for every $ \BW \in C^* \Lin $ of the form $ \BW = CF(\BX) $ for some $ \BX \in C^* \Lin $. 
In order to do this consider the set 
$$ 
C_\BW = [I, C^* \Lin]_{cof}((\BV_i)_{i \in I}, \Delta(\BW)) 
$$ 
of all cofibrant transformations between the $ I $-diagrams $ (\BV_i)_{i \in I} $ and $ \Delta(\BW) $. Note that $ C_\BW $ is nonempty since there 
exists a cofibrant transformation $ 0: (\BV_i)_{i \in I} \rightarrow \Delta(\BW) $, given by functors $ 0_i: \BV_i \rightarrow \BW $ 
mapping every $ V \in \BV_i $ to some zero object in $ \BW $, together with uniquely determined natural 
isomorphisms $ 0_{i,j}(i \to j): 0_i \Rightarrow 0_j \circ \iota(i \to j) $ for $ i \to j $ in $ I $. 

Let $ \kappa $ be a strong limit cardinal greater than the sum of the cardinalities of all categories $ \BV_i $ in the given family. 
We denote by $ C^* \Lin_\kappa $ the full sub $ 2 $-category of $ C^* \Lin $ given by a chosen set of $ C^* $-categories which contains all 
countably additive $ C^* $-categories of cardinality less than $ \kappa $ up to isomorphism and define 
$$
\BV = \prod_{\BW \in C^* \Lin_\kappa} \prod_{\sigma \in C_\BW} \BW.  
$$
Then $ \BV $ is a countably additive $ C^* $-category, 
and we obtain a canonical transformation $ t: (\BV_i)_{i \in I} \rightarrow \Delta(\BV) $ 
by stipulating that the composition of $ t $ with projection onto the component of $ \BV $ associated with $ \sigma \in C_\BW $ equals $ \sigma $. 
The transformation $ t $ consists of nondegenerate linear $ * $-functors $ t_i: \BV_i \rightarrow M\BV $ and unitary 
natural isomorphisms $ t_{i,j}(i \to j): t_i \Rightarrow t_j \circ \iota(i \to j) $, 
satisfying the required coherence conditions. 
All component transformations $ \sigma \in C_\BW $ are cofibrant by assumption, and therefore the transformation $ t $ is again cofibrant. 

Let $ \langle t \rangle \subset M\BV $ be the $ C^* $-category generated by the union of the images of all the functors $ t_i $ and the unitary 
natural isomorphisms $ t_{i,j}(i \to j) $. The objects of $ \langle t \rangle $ are all objects of the form $ t_i(V) $ for some $ i \in I $ 
and $ V \in \BV_i $, and the morphism spaces are generated by the images of all morphism spaces $ \BV_i(V, W) $ for $ V, W \in \BV_i $ under $ t_i $, as well 
as the multiplier morphisms $ t_{i,j}(i \to j)(V) $ for $ V \in \BV_i $ and their adjoints. That is, the morphism space $ \langle t \rangle(X, Y) $ 
is the norm closure of the linear span of all morphisms $ f_r \circ \cdots \circ f_1 \in M\BV(X,Y) $ consisting of composable strings of multiplier morphisms of the 
form $ f_k = t_{i_k,j_k}(i_k \to j_k)(V_k), f_k = t_{i_k,j_k}(i_k \to j_k)^*(V_k) $, or $ f_k = t_{i_k}(g_k) $ for some $ g_k \in \BV_{i_k}(V_k, W_k) $,  
with at least one morphism of the latter type. 
Using that $ \kappa $ is a strong limit cardinal one verifies that $ \langle t \rangle $ is a $ C^* $-category of cardinality less than $ \kappa $. 
By construction, the morphisms $ t_{i,j}(i \to j)(V) $ for $ V \in \BV_i $ 
are naturally contained in $ M \langle t \rangle(t_i(V), t_j \circ \iota(i \to j)(V)) $. 
Corestriction of $ t_i: \BV_i \rightarrow M\BV $ determines a nondegenerate linear $ * $-functor $ \BV_i \rightarrow \langle t \rangle $, which we will again 
denote by $ t_i $. 

Let $ \varinjlim_{i \in I} \BV_i = \ind(\langle t \rangle^{\oplus}) $ be the ind-completion of the finite additive completion of $ \langle t \rangle $, 
which can be viewed as the completion of $ \langle t \rangle $ under countable direct sums and subobjects. 
One checks that the category $ \varinjlim_{i \in I} \BV_i \in C^* \Lin $ has again cardinality less than $ \kappa $. 

The functors $ t_j $ induce nondegenerate linear $ * $-functors $ T_j: \BV_j \rightarrow \varinjlim_{i \in I} \BV_i $ for all $ j \in I $. 
Moreover the unitary natural transformations $ t_{i,j}(i \to j) $ induce unitary 
natural transformations $ T_{i,j}(i \to j): T_i \Rightarrow T_j \circ \iota(i \to j) $, 
assembling to a cofibrant transformation $ T: (\BV_i)_{i \in I} \rightarrow \Delta(\varinjlim_{i \in I} \BV_i) $. 

Let us next show that precomposition with $ T $ induces an equivalence of $ C^* $-categories 
$$
F_\BW: C^* \Lin(\varinjlim_{i \in I} \BV_i, \BW) \rightarrow [I, C^* \Lin]_{cof}((\BV_i)_{i \in I}, \Delta(\BW)) 
$$
for every $ \BW \in C^* \Lin $ of cardinality less than $ \kappa $. For this it suffices to consider $ \BW \in C^* \Lin_\kappa $.  

Assume that $ \sigma: (\BV_i)_{i \in I} \rightarrow \Delta(\BW) $ is a transformation. By the construction of $ \BV $, projection onto the component 
corresponding to $ \sigma $ in $ \BV $ induces a nondegenerate linear $ * $-functor $ \Sigma: \varinjlim_{i \in I} \BV_i \rightarrow \BW $, 
such that $ \Delta(\Sigma) \circ T = \sigma $. It follows that $ F_\BW $ is essentially surjective on objects. 

Next assume that $ \Bf, \Bg: \varinjlim_{i \in I} \BV_i \rightarrow M\BW $ are nondegenerate linear $ * $-functors, 
and let $ \phi: \Bf \Rightarrow \Bg $ be a multiplier natural transformation such that $ F_\BW(\phi) = \Delta(\phi) * T = 0 $. 
Then $ \phi(T_j(V_j)): \Bf(T_j(V_j)) \rightarrow \Bg(T_j(V_j)) $ 
vanishes for all $ j \in I $ and $ V_j \in \BV_j $. Every object of $ \varinjlim_{i \in I} \BV_i $ is a subobject of a countable direct sum of such 
objects $ T_j(V_j) $. Since both $ \Bf, \Bg $ are nondegenerate linear $ * $-functors it follows from naturality that $ \phi(V) = 0 $ for 
all $ V \in \varinjlim_{i \in I} \BV_i $. 
Hence $ F_\BW $ is faithful. 

Now assume that $ \Gamma: F_\BW(\Bf) \Rightarrow F_\BW(\Bg) $ is a modification.  
We claim that we can assemble the multiplier morphisms $ \Gamma_i(V_i): \Bf(T_i(V_i)) \Rightarrow \Bg(T_i(V_i)) $ for $ V_i \in \BV_i $ 
to a multiplier natural transformation $ \Phi: \Bf \Rightarrow \Bg $ such that $ F_\BW(\Phi) = \Gamma $. 
Since the union of the objects in $ \BV_i $ for $ i \in I $ gets mapped injectively into $ \BV $ under the functors $ T_i $, 
the morphisms $ \Gamma_i(V_i) $ assemble uniquely to define a uniformly bounded family of multiplier morphisms 
$ \Phi(T_i(V_i)) = \Gamma(V_i): \Bf(T_i(V_i)) \rightarrow \Bg(T_i(V_i)) $ for $ V_i \in \BV_i $, 
and this extends canonically to direct sums and their subobjects. 
Using the fact that $ \Gamma $ is a modification one checks that the resulting uniformly bounded family of multiplier morphisms 
$ \Phi(V): \Bf(V) \rightarrow \Bg(V) $ for $ V \in \varinjlim_i \BV_i $ defines a multiplier natural transformation $ \Phi: \Bf \Rightarrow \Bg $ as required. 
It follows that $ F_\BW $ is full. 

Suppose that we replace $ \kappa $ in the above constructions by some strong limit cardinal $ \lambda \geq \kappa $. Direct inspection shows 
the resulting object $ \varinjlim^\lambda_{i \in I} \BV_i $ in $ C^* \Lin $ has again cardinality strictly less than $ \kappa $, observing that this 
cardinality is determined entirely in terms of the $ I $-diagram $ (\BV_i)_{i \in I} $. 
Therefore $ \varinjlim^\lambda_{i \in I} \BV_i $ is isomorphic to an object of $ C^* \Lin_\kappa $, and our above 
considerations imply that it must be equivalent to $ \varinjlim_{i \in I} \BV_i $. 

Finally, let $ \BW \in C^* \Lin $ be an arbitrary countably additive $ C^* $-category. Then $ \BW $ is isomorphic to an object of $ C^* \Lin_\lambda $ for 
a sufficiently large strong limit cardinal $ \lambda \geq \kappa $. By the above reasoning we see that $ F_\BW $ induces an equivalence 
$$
C^* \Lin(\varinjlim_{i \in I} \BV_i, \BW) \simeq [I, C^* \Lin]_{cof}((\BV_i)_{i \in I}, \Delta(\BW)) 
$$
as desired. \qed 

As a special case of Theorem \ref{bicolimitexistence} we see that $ C^* \Lin $ admits all bicoproducts, bicoequalisers and bipushouts. 

Bicoproducts can be described concretely as follows. Let $ I $ be an index set and let $ (\BV_i)_{i \in I} $ be a family of countably additive $ C^* $-categories. 
The \emph{direct sum} $ \bigoplus_{i \in I} \BV_i $ is the full subcategory of $ \prod_{i \in I} \BV_i $ whose objects are families $ (X_i)_{i \in I} $ 
of objects $ X_i \in \BV_i $, with at most countably many $ X_i $ being nonzero. This is again a countably additive $ C^* $-category. 

We obtain canonical fully faithful linear $ * $-functors $ \iota_j: \BV_j \rightarrow \bigoplus_{i \in I} \BV_i $ for all $ j \in I $, 
mapping $ V \in \BV_j $ to the family whose only nonzero object is $ V $, based at $ j \in I $. 

\begin{prop} \label{categoricalcstardirectsums}
Let $ I $ be an index set and $ (\BV_i)_{i \in I} $ a family of $ C^* $-categories in $ C^* \Lin $. 
For every $ \BW \in C^* \Lin $, precomposition with the family of linear $ * $-functors $ \iota_j $ induces an equivalence of $ C^* $-categories 
$$
C^* \Lin(\bigoplus_{i \in I} \BV_i, \BW) \simeq M \prod_{i \in I} C^* \Lin(\BV_i, \BW) = [I, C^* \Lin]((\BV_i)_{i \in I}, \Delta(\BW)), 
$$
pseudonatural in $ \BW $. 
\end{prop} 

\proof The assignment $ F: C^* \Lin(\bigoplus_{i \in I} \BV_i, \BW) \rightarrow M\prod_{i \in I} C^* \Lin(\BV_i, \BW) $ given by 
$ F(\Bf) = (\Bf \circ \iota_i)_{i \in I} $ on objects and $ F(\sigma) = (\sigma * \iota_i)_{i \in I} $ on morphisms yields a well-defined unital linear $ * $-functor. 

If $ (\Bh_i)_{i \in I} $ is a family of nondegenerate linear $ * $-functors $ \Bh_i: \BV_i \rightarrow \BW $ we obtain a nondegenerate linear 
$ * $-functor $ \Bh: \bigoplus_{i \in I} \BV_i \rightarrow \BW $ by setting $ \Bh((V_i)_{i \in I}) = \bigoplus_{i \in I} \Bh_i(V_i) $ on objects
and $ \Bh((f_i)_{i \in I}) = \bigoplus_{i \in I} \Bh_i(f_i) $ on morphisms. Here we use that $ \BW $ has countable direct sums and 
that $ (V_i)_{i \in I} \in \bigoplus_{i \in I} \BV_i $ has at most countably many nonzero entries. 
We have $ F(\Bh) \cong (\Bh_i)_{i \in I} $ by construction, and it follows that $ F $ is essentially surjective. 

If $ (\sigma_i)_{i \in I} $ is a uniformly bounded family of multiplier natural transformations $ \sigma_i: \Bh \circ \iota_i \Rightarrow \Bk \circ \iota_i $ 
then $ \sigma((V_i)_{i \in I}) = \prod_{i \in I} \sigma_i(V_i) $ defines a multiplier natural transformation $ \sigma: \Bh \Rightarrow \Bk $ such 
that $ F(\sigma) = (\sigma_i)_{i \in I} $. 
It follows that $ F $ is full. Finally, note that $ F(\sigma) = 0 $ for $ \sigma: \Bh \Rightarrow \Bk $ means $ \sigma(\iota_i(V_i)) = 0 $ for 
all $ i \in I $ and $ V_i \in \BV_i $, and therefore $ \sigma(V) = 0 $ for all $ V \in \bigoplus_{i \in I} \BV_i $ by naturality. Hence $ F $ is faithful. \qed 

Let us briefly explain the connection between Theorem \ref{bicolimitexistence} and the results obtained in \cite{AlbandikMeyercolimits}. 
We shall not aim for greatest generality here in order not to obscure the argument.  

Let us consider a diagram $ \iota: I \rightarrow C^* \Lin $ such that the index $ 1 $-category $ I $ is countable, we have $ \iota(i) = \BV_i = \Hilb_{A_i} $ 
for separable $ C^* $-algebras $ A_i $ for all $ i \in I $, and all nondegenerate linear $ * $-functors $ \iota(i \to j): \BV_i \rightarrow M\BV_j $ 
for $ i \to j $ in $ I $ are given by proper $ A_i $-$ A_j $-correspondences in $ \Hilb_{A_j} $. 
Then we can view $ \iota $ as a diagram both in $ C^* \Lin $ and in the (proper) correspondence bicategory.  
Let $ \O $ be the $ C^* $-algebra constructed for this diagram in \cite{AlbandikMeyercolimits}. 

\begin{prop} \label{AMcomparison}
Under the above hypotheses we have an equivalence
$$ 
\varinjlim_{i \in I} \BV_i \simeq \Hilb_\O 
$$
of $ C^* $-categories. 
\end{prop} 

\proof Under our assumptions, an inspection of the constructions in \cite{AlbandikMeyercolimits} shows that $ \O $ is separable, 
and that the defining $ A_i $-$ \O $-correspondences are countably generated and proper.
On the other hand, since $ I $ is countable it follows from the constructions in the proof of Theorem \ref{bicolimitexistence} that $ \varinjlim_{i \in I} \BV_i $ 
has separable morphism spaces and is singly generated. 
If we write $ \varinjlim_{i \in I} \BV_i \simeq \Hilb_A $ for a generator $ A $, then 
the $ C^* $-algebras $ A $ and $ \O $ satisfy the same universal property in the bicategory of separable $ C^* $-algebras and countably generated correspondences. 
As such they have to equivalent, which yields the claim. \qed 

In other words, Proposition \ref{AMcomparison} shows that the $ C^* $-algebra $ \O $ provides a concrete model for the bicolimit of the 
diagram $ \iota: I \rightarrow C^* \Lin $. Together with the results in \cite{AlbandikMeyercolimits}, this yields a rich supply of concrete examples of bicolimits 
of linear $ C^* $-categories. 

Let us finish this section by noting that the $ 2 $-category $ C^* \lin $ of finitely additive $ C^* $-categories admits bicolimits as well. Since 
the proof of this assertion is largely parallel to the proof of Theorem \ref{bicolimitexistence} we shall only state the result as follows. 

\begin{theorem} 
Let $ I $ be a small category and let $ (\BV_i)_{i \in I} $ be an $ I $-diagram in $ C^* \lin $. Then there exists a 
bicolimit $ \varinjlim_{i \in I} \BV_i \in C^* \lin $, that is, there exists a transformation $ T: (\BV_i)_{i \in I} \rightarrow \Delta(\varinjlim_{i \in I} \BV_i) $
such that precomposition with $ T $ induces an equivalence of $ C^* $-categories
$$
C^* \lin(\varinjlim_{i \in I} \BV_i, \BW) \simeq [I, C^* \lin]((\BV_i)_{i \in I}, \Delta(\BW))
$$
for every $ \BW \in C^* \lin $, pseudonatural in $ \BW $. 
\end{theorem}

\section{Balanced tensor products} \label{secbalancedtensorproducts}

In this section we discuss balanced tensor products of module categories over $ C^* $-tensor categories. 

Let us first specify our setup. 

\begin{definition} \label{deflinearcstartensor}
A countably additive $ C^* $-tensor category is a $ C^* $-category $ \BA \in C^* \Lin $ together with 
\begin{bnum}
\item[$ \bullet $] a nondegenerate bilinear $ * $-functor $ \otimes: \BA \times \BA \rightarrow M\BA $, 
\item[$ \bullet $] an object $ 1 \in \BA $, 
\item[$ \bullet $] a unitary 
natural isomorphism 
$$ 
\alpha: \otimes \circ (\otimes \times \id) \Rightarrow \otimes \circ (\id \times \otimes) 
$$ 
written 
$$ 
\alpha_{X, Y, Z}: (X \otimes Y) \otimes Z \rightarrow X \otimes (Y \otimes Z)
$$ 
for all $ X, Y, Z \in \BA $, and called associator, 
\item[$ \bullet $] unitary natural isomorphisms $ \rho: - \otimes 1 \Rightarrow \id $ 
and $ \lambda: 1 \otimes - \Rightarrow \id $, called left and right unitors, written $ \rho_X: X \otimes 1 \rightarrow X $ 
and $ \lambda_X: 1 \otimes X \rightarrow X $ for $ X \in \BA $, respectively. 
\end{bnum}
These data are supposed to satisfy the following conditions. 
\begin{bnum}
\item[$\bullet$] (Associativity constraints) For all objects $ W, X, Y, Z \in \BA $ the diagram
\begin{center}
\begin{tikzcd}[column sep = -12mm]
			&&
			(W \otimes X) \otimes (Y \otimes Z)
			\arrow[drr,"\alpha_{W, X, Y \otimes Z}"]
			&&
			\\
			((W \otimes X) \otimes Y) \otimes Z
			\arrow[urr,"\alpha_{W \otimes X, Y, Z}"]
			\arrow[dr, swap, "\alpha_{W, X, Y} \otimes \id"]
			&&&&
			W \otimes (X \otimes (Y \otimes Z))
			\\
			&
			(W \otimes (X \otimes Y)) \otimes Z
			\arrow[rr, swap, "\alpha_{W, X \otimes Y, Z}"]
			&&
			W \otimes ((X \otimes Y) \otimes Z)
			\arrow[ur, swap, "\id \otimes \alpha_{X, Y, Z}"]
			&
\end{tikzcd}
\end{center}
is commutative. 
\item[$\bullet$] (Unit constraints) $ \lambda_1 = \rho_1 $ and for all objects $ X, Y \in \BA $ the diagram 
\begin{center}
\begin{tikzcd}
			(X \otimes 1) \otimes Y
			\arrow[rr, "\alpha_{X, 1, Y}"]
			\arrow[dr, swap, "\rho_X \otimes \id"]
			&&
			X \otimes (1 \otimes Y)
			\arrow[dl, "\id \otimes \lambda_Y"]
			\\
			&
			X \otimes Y
\end{tikzcd}
\end{center}
is commutative.  
\end{bnum} 
\end{definition} 

Definition \ref{deflinearcstartensor} means that countably additive $ C^* $-tensor categories are monoids 
in the monoidal $ 2 $-category $ C^* \Lin $ of countably additive $ C^* $-categories with the maximal tensor product, compare Theorem \ref{symmon2cat}. 
We will also consider $ C^* \lin $, leading to the following variant of Definition \ref{deflinearcstartensor}. 

\begin{definition} \label{deffinitelylinearcstartensor}
A finitely additive $ C^* $-tensor category is a category $ \BA \in C^* \lin $ together with  
\begin{bnum}
\item[$ \bullet $] a unital bilinear $ * $-functor $ \otimes: \BA \times \BA \rightarrow \BA $, 
\item[$ \bullet $] an object $ 1 \in \BA $, 
\item[$ \bullet $] a unitary natural isomorphism 
$ \alpha: \otimes \circ (\otimes \times \id) \Rightarrow \otimes \circ (\id \times \otimes) $ 
\item[$ \bullet $] unitary natural isomorphisms $ \rho: - \otimes 1 \Rightarrow \id $ and $ \lambda: 1 \otimes - \Rightarrow \id $, 
\end{bnum} 
satisfying the same associativity and unit constraints as in the countably additive case.  
\end{definition} 

Definition \ref{deffinitelylinearcstartensor} agrees with the notion of a $ C^* $-tensor category used in \cite{NTlecturenotes} 
except that we do not require unit objects to be simple. 

\begin{example} \label{finitelylinearexamples}
Let us list some examples of finitely additive $ C^* $-tensor categories. 
\begin{bnum} 
\item[a)] If $ G $ is a compact (quantum) group then the representation category $ \Rep(G) $ of $ G $ is a finitely additive $ C^* $-tensor category. 
The objects of $ \Rep(G) $ are finite dimensional unitary $ G $-representations, and morphisms are all $ G $-equivariant linear operators. 
The tensor structure is given by the tensor product of $ G $-representations, which for classical groups is the usual diagonal action on the tensor 
product of the underlying Hilbert spaces. Note that all morphism spaces in $ \Rep(G) $ are finite dimensional, 
and that the tensor unit is given by the trivial representation on $ \mathbb{C} $. In particular $ \Rep(G) $ is semisimple with simple tensor unit. \\
More generally one obtains rigid $ C^* $-tensor categories from bimodule categories of finite index subfactors, 
or from Hilbert $ C^* $-bimodules in the sense of \cite{KajiwaraWatatanihilbertsctarbimodules}. 
In particular, the category of $ M $-bimodules generated by a finite index subfactor $ N \subset M $ is a semisimple finitely additive $ C^* $-tensor category
with simple tensor unit. 

\item[b)] Let $ A $ be a unital $ C^* $-algebra and consider the category $ _A \Hilb^f_A $ of finite right Hilbert $ A $-bimodules. 
By definition, the objects of $ _A \Hilb^f_A $ are finitely generated projective right Hilbert $ A $-modules $ \E \in \Hilb^f_A $ 
equipped with a unital $ * $-representation $ A \rightarrow \LH(\E) = \KH(\E) $. Morphisms in $ _A \Hilb^f_A $ are adjointable operators 
which commute with the left actions. 
The tensor structure on $ _A \Hilb^f_A $ is given by the balanced tensor product of Hilbert modules, and the tensor unit is $ A \in \Hilb^f_A $ viewed as 
a bimodule with the left multiplication action. 
The resulting finitely additive $ C^* $-tensor category is not semisimple in general. 
Morphism spaces in $ _A \Hilb^f_A $ are typically infinite dimensional, and the tensor unit is not simple. 

\item[c)] Let $ \BV \in C^* \lin $ be a finitely additive $ C^* $-category. Then the category $ C^* \lin(\BV, \BV) $ of all unital 
linear $ * $-functors $ \BV \rightarrow \BV $ is a finitely additive $ C^* $-tensor category with the tensor structure induced by composition of functors. 
In fact, up to equivalence the category $ _A \Hilb^f_A $ described in the previous example is the special case of this construction for $ \BV = \Hilb^f_A $. 
\end{bnum} 
\end{example}

Using ind-completion one can pass from finitely additive $ C^* $-tensor categories to countably additive $ C^* $-tensor categories as follows. 

\begin{lemma} \label{finitelylineartolinear} 
Let $ \BA $ be a finitely additive $ C^* $-tensor category. Then $ \ind \BA $ is naturally a countably additive $ C^* $-tensor category 
such that the canonical inclusion functor $ \BA \rightarrow \ind \BA $ is compatible with tensor products. 
\end{lemma} 

\proof The tensor structure for $ \BA $ can be viewed as a unital linear $ * $-functor $ \otimes: \BA \boxtimes_{\max} \BA \rightarrow \BA $. 
It extends canonically to a nondegenerate linear $ * $-functor $ \ind(\BA \boxtimes_{\max} \BA) \rightarrow \ind \BA $ on the level of ind-categories. 
Combining this with the equivalence $ \ind \BA \boxtimes_{\max} \ind \BA \simeq \ind(\BA \boxtimes_{\max} \BA) $ from Proposition \ref{finitetensorproducts}  
gives a nondegenerate linear $ * $-functor $ \ind \BA \boxtimes_{\max} \ind \BA \rightarrow \ind \BA $, which we denote again by $ \otimes $. 
The tensor unit $ 1 \in \BA \subset \ind \BA $ and the extensions of $ \alpha, \lambda, \rho $ to the ind-categories provide the data 
of a countably additive $ C^* $-tensor category structure on $ \ind \BA $, and one checks that the associativity and unit constraints continue to hold. \qed 

Applying Lemma \ref{finitelylineartolinear} to the examples in \ref{finitelylinearexamples} yields basic examples of countably additive $ C^* $-tensor categories. 
The ind-categories of representation categories of compact quantum groups and bimodule categories of subfactors feature naturally in the study of approximation 
properties \cite{NYdrinfeldcenter}. 

\begin{example} \label{linearexamples}
Let us describe some examples of countably additive $ C^* $-tensor categories which are not obtained as ind-categories of finitely additive $ C^* $-tensor categories. 
\begin{bnum} 
\item[a)] Let $ X $ be a locally compact Hausdorff space and let $ A = C_0(X) $. The category of symmetric right Hilbert $ A $-bimodules 
has as objects Hilbert modules $ \E \in \Hilb_A $, viewed as symmetric $ A $-bimodules via $ a \cdot \xi = \xi a $ for $ a \in A, \xi \in \E $, 
and morphisms as in $ \Hilb_A $. 
The interior tensor product of Hilbert modules turns this into a countably additive $ C^* $-tensor category with unit object $ A $. 
If $ X $ is noncompact this category will typically fail to be finitely presented. 

\item[b)] Let $ X $ be a locally compact Hausdorff space, let $ A_0 = C_0(X) $, and let $ A = A_1 $ be a continuous trace algebra with spectrum $ X $. 
For $ n > 1 $ let 
$$ 
A_n = A^{\otimes_{C_0(X)} n} = A \otimes_{C_0(X)} A \otimes_{C_0(X)} \cdots \otimes_{C_0(X)} A 
$$ 
be the the $ n $-fold (maximal) $ C_0(X) $-tensor product, compare \cite{Blancharddef}. Then 
$$ 
\BA = \bigoplus_{n = 0}^\infty \Hilb_{A_n} 
$$ 
is a $ C^* $-tensor category with the tensor product given by $ \E \otimes_{C_0(X)} \F \in \Hilb_{A_{m + n}} $ for $ \E \in \Hilb_{A_m}, \F \in \Hilb_{A_n} $. 
This is well-defined and associative up to isomorphism. The tensor unit is $ C_0(X) $ viewed as Hilbert module over itself. 

\item[c)] Let $ \BA $ be a countably additive $ C^* $-tensor category and consider 
$$
\KH_\BA = \bigoplus_{i,j = 1}^\infty \BA 
$$
as a $ C^* $-tensor category with the categorical version of matrix multiplication. 
That is, denoting $ \BA_{ij} \subset \KH_\BA $ the copy of $ \BA $ corresponding to the indices $ i,j \in \mathbb{N} $, the tensor 
product functor $ \KH_\BA \times \KH_\BA \rightarrow \KH_\BA $ maps $ \BA_{ij} \times \BA_{kl} $ 
to zero if $ j \neq k $, and by the given tensor structure $ \BA \times \BA \rightarrow M(\BA) $ to $ M\BA_{il} \subset M(\KH_\BA) $ if $ j = k $. 
The unit object of $ \KH_\BA $ has the tensor unit $ 1 \in \BA $ on the diagonal, and the zero objects else. 
Moreover, the associativity and unit constraints of $ \BA $ induces associativity and unit constraints for $ \KH_\BA $ in a natural way. 
We shall refer to $ \KH_\BA $ as the infinite matrix $ C^* $-tensor category over $ \BA $. If $ \BA $ is not finitely presented then neither is $ \KH_\BA $. 
\end{bnum} 
\end{example}

Let us next discuss module categories, module functors and their natural transformations, compare \cite{Ostrikmodulecategories}. Since the definitions 
are parallel for finitely and countably additive $ C^* $-categories, we will only state the countably additive case. 

\begin{definition} 
Let $ \BA $ be a countably additive $ C^* $-tensor category. A left $ \BA $-module category is a $ C^* $-category $ \BV \in C^* \Lin $ together with 
\begin{bnum}
\item[$ \bullet $] a nondegenerate bilinear $ * $-functor $ \otimes: \BA \times \BV \rightarrow M\BV $ 
\item[$ \bullet $] a unitary multiplier natural isomorphism 
$$ 
\alpha: \otimes \circ (\otimes \times \id) \Rightarrow \otimes \circ (\id \times \otimes)
$$ 
written 
$$ 
\alpha_{X, Y, V}: (X \otimes Y) \otimes V \rightarrow X \otimes (Y \otimes V)
$$ 
for all $ X, Y \in \BA, V \in \BV $, and called associator. 
\item[$ \bullet $] a unitary multiplier natural isomorphism $ \lambda: 1 \otimes - \Rightarrow \id $, called unitor, 
and written $ \lambda_V: 1 \otimes V \rightarrow V $ for $ V \in \BV $. 
\end{bnum}
These data are supposed to satisfy the following conditions. 
\begin{bnum}
\item[$\bullet$] (Associativity constraints) For all objects $ W, X, Y, \in \BA, V \in \BV $ the diagram
\begin{center}
\begin{tikzcd}[column sep = -12mm]
			&&
			(W \otimes X) \otimes (Y \otimes V)
			\arrow[drr,"\alpha_{W, X, Y \otimes V}"]
			&&
			\\
			((W \otimes X) \otimes Y) \otimes V
			\arrow[urr,"\alpha_{W \otimes X, Y, V}"]
			\arrow[dr, swap, "\alpha_{W, X, Y} \otimes \id"]
			&&&&
			W \otimes (X \otimes (Y \otimes V))
			\\
			&
			(W \otimes (X \otimes Y)) \otimes V
			\arrow[rr, swap, "\alpha_{W, X \otimes Y, V}"]
			&&
			W \otimes ((X \otimes Y) \otimes V)
			\arrow[ur, swap, "\id \otimes \alpha_{X, Y, V}"]
			&
\end{tikzcd}
\end{center}
is commutative. 
\item[$\bullet$] (Unit constraints) For all objects $ X \in \BA, V \in \BV $ the diagram 
\begin{center}
\begin{tikzcd}
			(X \otimes 1) \otimes V
			\arrow[rr, "\alpha_{X, 1, V}"]
			\arrow[dr, swap, "\rho_X \otimes \id"]
			&&
			X \otimes (1 \otimes V)
			\arrow[dl, "\id \otimes \lambda_V"]
			\\
			&
			X \otimes V
\end{tikzcd}
\end{center}
is commutative.  
\end{bnum} 
\end{definition} 

A module functor between $ \BA $-module categories $ \BV, \BW $ is a nondegenerate linear $ * $-functor $ \Bf: \BV \rightarrow M\BW $ together with a 
unitary multiplier natural isomorphism $ \gamma: \Bf \circ \otimes \Rightarrow \otimes \circ (\id \times \Bf) $, 
written $ \gamma_{X, V}: \Bf(X \otimes V) \rightarrow X \otimes \Bf(V) $, such that the diagrams 
\begin{center}
\begin{tikzcd}[column sep = -5mm]
			&&
			(X \otimes Y) \otimes \Bf(V)
			\arrow[drr,"\alpha_{X, Y, \Bf(V)}"]
			&&
			\\
			\Bf((X \otimes X) \otimes V) 
			\arrow[urr,"\gamma_{X \otimes Y, V}"]
			\arrow[dr, swap, "\Bf(\alpha_{X, Y, V})"]
			&&&&
			X \otimes (Y \otimes \Bf(V))
			\\
			&
			\Bf(X \otimes (Y \otimes V))
			\arrow[rr, swap, "\gamma_{X, Y \otimes V}"]
			&&
			X \otimes \Bf(Y \otimes V)
			\arrow[ur, swap, "\id \otimes \gamma_{Y, V}"]
			&
\end{tikzcd}
\end{center}
and 
\begin{center}
\begin{tikzcd}
			\Bf(1 \otimes V)
			\arrow[rr, "\gamma_{1, V}"]
			\arrow[dr, swap, "\Bf(\lambda_V)"]
			&&
			1 \otimes \Bf(V)
			\arrow[dl, "\lambda_{\Bf(V)}"]
			\\
			&
			\Bf(V)
\end{tikzcd}
\end{center}
are commutative. 

An $ \BA $-module natural transformation between $ \BA $-module functors $ \Bf, \Bg: \BV \rightarrow M\BW $, with 
corresponding $ \gamma: \Bf \circ \otimes \Rightarrow \otimes \circ (\id \times \Bf), \eta: \Bg \circ \otimes \Rightarrow \otimes \circ (\id \times \Bg) $, 
is a multiplier natural transformation $ \phi: \Bf \Rightarrow \Bg $ of the underlying nondegenerate linear $ * $-functors such that 
\begin{center}
\begin{tikzcd}
			\Bf(X \otimes V)
		  \arrow[r, "\gamma_{X, V}"]
			\arrow[d, swap, "\phi_{X \otimes V}"]
			&
			X \otimes \Bf(V)
			\arrow[d, "\id \otimes \phi_V"]
			\\
			\Bg(X \otimes V)
			\arrow[r, "\eta_{X, V}"]
			&
			X \otimes \Bg(V)
\end{tikzcd}
\end{center}
is commutative for all $ X \in \BA $ and $ V \in \BV $. 

In a similar way one defines right $ \BA $-module categories, their module functors, and module natural transformations. 
More generally, if $ \BA, \BB $ are countably additive $ C^* $-tensor categories then one may define $ \BA $-$ \BB $-bimodule categories, bimodule functors 
and bimodule natural transformations, compare for instance \cite{ENOfusionhomotopy}, \cite{Greenoughmon} in the algebraic situation. 
Since bimodule categories are not strictly needed for our purposes we shall not write down these definitions. 
 
\begin{example} \label{moduleexamples}
Let us give some examples of module categories and bimodule categories. 
\begin{bnum} 
\item[a)] Every countably additive $ C^* $-tensor category is a left and right module category over itself via the tensor structure, 
and in fact a bimodule category. 
\item[b)] Let $ A $ be a unital $ C^* $-algebra and let $ _A\Hilb^f_A $ be the finitely additive $ C^* $-tensor category from Example \ref{finitelylinearexamples} b). 
If $ B $ is a unital $ C^* $-algebra then the category $ _A \Hilb^f_B $ of finitely generated projective right Hilbert $ B $-modules $ \E $ 
equipped with a unital $ * $-representation $ A \rightarrow \LH(\E) $ and $ A $-linear compact operators is a finitely additive left module 
category over $ _A\Hilb^f_A $ with the action given by interior tensor product. 
\item[c)] Let $ \BA $ be a countably additive $ C^* $-tensor category, and let $ \KH_\BA $ the countably additive $ C^* $-tensor category from 
Example \ref{linearexamples} c). Applying the categorical version of matrix multiplication we obtain a $ \KH_\BA $-$ \BA $-bimodule category 
$$ 
\HH_\BA = \bigoplus_{n = 1}^\infty \BA, 
$$
viewing the objects of $ \HH_\BA $ as column vectors of objects from $ \BA $. 
Dually, by taking row vectors instead we obtain a $ \BA $-$ \KH_\BA $-bimodule category $ \HH_\BA^\vee = \bigoplus_{n = 1}^\infty \BA $ 
with the same underlying $ C^* $-category. 
\end{bnum} 
\end{example} 

Let us now discuss balanced bilinear functors and balanced tensor products. As above, we only state the case of countably additive $ C^* $-categories, 
since the finitely additive case is analogous. 

\begin{definition} \label{defbalancedfunctor}
Let $ \BA $ be a countably additive $ C^* $-tensor category. Moreover let $ \BV $ be a right $ \BA $-module category, $ \BW $ a left $ \BA $-module category, 
and $ \BX \in C^* \Lin $. An $ \BA $-balanced functor $ \Bf: \BV \times \BW \rightarrow M\BX $ is a nondegenerate bilinear $ * $-functor 
together with a unitary multiplier natural isomorphism $ \beta: \Bf \circ (\otimes \times \id_\BW) \Rightarrow \Bf \circ (\id_\BV \times \otimes) $,  
written $ \beta_{V, A, W}: \Bf(V \otimes A, W) \rightarrow \Bf(V, A \otimes W) $, such that 
for all $ V \in \BV, A, B \in \BA, W \in \BW $ the diagram 
\begin{center}
\begin{tikzcd}[column sep = -10mm]
			&&
			\Bf(V \otimes A, B \otimes W)
			\arrow[drr,"\beta_{V, A, B \otimes W}"]
			&&
			\\
			\Bf((V \otimes A) \otimes B, W) 
			\arrow[urr,"\beta_{V \otimes A, B, W}"]
			\arrow[dr, swap, "\id_\Bf * (\alpha_{V, A, B} \times \id)"]
			&&&&
			\Bf(V, A \otimes (B \otimes W))
			\\
			&
			\Bf(V \otimes (A \otimes B), W)
			\arrow[rr, swap, "\beta_{V, A \otimes B, W}"]
			&&
			\Bf(V, (A \otimes B) \otimes W)
			\arrow[ur, swap, "\id_\Bf * (\id \times \alpha_{A, B, W})"]
			&
\end{tikzcd}
\end{center}
is commutative. 
\end{definition} 

Sometimes $ \BA $-balanced functors are also called $ \BA $-bilinear.
If $ \BA = \Hilb $ is the $ C^* $-tensor category of separable Hilbert spaces then any bilinear $ * $-functor is automatically $ \BA $-balanced. 
According to Proposition \ref{universalpropmaxtensor} every $ \BA $-balanced functor $ \BV \times \BW \rightarrow M\BX $ factorises 
through $ \BV \boxtimes_{\max} \BW $, since being $ \BA $-balanced is further structure and properties on top of a nondegenerate bilinear $ * $-functor. 
In the sequel we will often identify $ \BA $-balanced functors $ \BV \times \BW \rightarrow M\BX $ with nondegenerate 
linear $ * $-functors $ \BV \boxtimes_{\max} \BW \rightarrow M\BX $. 

\begin{definition} \label{defbalancedtrans}
Let $ \Bf, \Bg: \BV \times \BW \rightarrow \BX $ be $ \BA $-balanced bilinear $ * $-functors, with balancing transformations 
$ \beta: \Bf \circ (\otimes \times \id_\BW) \Rightarrow \Bf \circ (\id_\BV \times \otimes), 
\gamma: \Bg \circ (\otimes \times \id_\BW) \Rightarrow \Bg \circ (\id_\BV \times \otimes) $, respectively.
A multiplier natural transformation $ \phi: \Bf \Rightarrow \Bg $ is called $ \BA $-balanced if 
\begin{center}
\begin{tikzcd}
			\Bf(V \otimes A, W)
			\arrow[r, "\beta_{V,A,W}"]
			\arrow[d, swap, "\phi_{V \otimes A, W}"]
			&
			\Bf(V, A \otimes W)
			\arrow[d, "\phi_{V, A \otimes W}"]
			\\
			\Bg(V \otimes A, W)
			\arrow[r, "\gamma_{V,A,W}"]
			&
			\Bg(V, A \otimes W)
\end{tikzcd}
\end{center}
is commutative for all $ V \in \BV, W \in \BW, A \in \BA $. 
\end{definition} 

If $ \BV $ is a right $ \BA $-module category, $ \BW $ a left $ \BA $-module category and $ \BX \in C^* \Lin $, then $ \BA $-balanced 
functors $ \BV \times \BW \rightarrow \BA $ together with the $ \BA $-balanced multiplier natural transformations form a $ C^* $-category 
which we denote by $ C^* \Bilin_\BA(\BV, \BW; \BX) $. 

Let us now state the main result of this section. 

\begin{theorem} \label{balancedcstartensor}
Let $ \BA \in C^* \Lin $ be a countably additive $ C^* $-tensor category. Morever let $ \BV_1, \BV_2 \in C^* \Lin $ be right and left $ \BA $-module categories, respectively. Then there exists a $ C^* $-category $ \BV_1 \boxtimes_\BA \BV_2 \in C^* \Lin $ together with an $ \BA $-balanced 
bilinear $ * $-functor $ \boxtimes_\BA: \BV_1 \times \BV_2 \rightarrow \BV_1 \boxtimes_\BA \BV_2 $ such that precomposition with $ \boxtimes_\BA $ 
induces an equivalence of $ C^* $-categories
$$
C^* \Bilin_\BA(\BV_1, \BV_2; \BW) \simeq C^* \Lin(\BV_1 \boxtimes_\BA \BV_2, \BW)
$$
for all $ \BW \in C^* \Lin $, pseudonatural in $ \BW $. 
\end{theorem} 

\proof Throughout the proof we shall abbreviate $ \boxtimes_{\max} = \boxtimes $. Moreover we will write $ \otimes: \BA \boxtimes \BA \rightarrow M\BA $ 
for the tensor product functor of $ \BA $ and denote by $ \lambda: \BA \boxtimes \BV_2 \rightarrow M\BV_2 $ and $ \rho: \BV_1 \boxtimes \BA \rightarrow M\BV_1 $ 
the module category actions. 

Let $ I $ be the opposite of the $ 2 $-truncated presimplicial category. This means that $ I $ is the category with three 
objects $ \bra 2 \ket, \bra 1 \ket, \bra 0 \ket $ 
and morphism sets generated by morphisms $ \partial_i: \bra n \ket \rightarrow \bra n - 1 \ket $ for $ 0 \leq i \leq n $ such 
that $ \partial_i \circ \partial_j = \partial_{j - 1} \circ \partial_i $ for $ i < j $. 

We consider the $ I $-diagram $ \iota: I \rightarrow C^* \Lin $ obtained from the truncated Bar complex 
\begin{center}
\begin{tikzcd}
			\BV_1 \boxtimes \BA \boxtimes \BA \boxtimes \BV_2
			\arrow[r, "", start anchor = {[yshift = 2mm]}, end anchor = {[yshift = 2mm]}]
			\arrow[r, ""]
    	\arrow[r, swap, "", start anchor = {[yshift = -2mm]}, end anchor = {[yshift = -2mm]}]
			&
			\BV_1 \boxtimes \BA \boxtimes \BV_2
			\arrow[r, "", start anchor = {[yshift = 1mm]}, end anchor = {[yshift = 1mm]}]
    	\arrow[r, swap, "", start anchor = {[yshift = -1mm]}, end anchor = {[yshift = -1mm]}]
			&
			\BV_1 \boxtimes \BV_2 
\end{tikzcd}
\end{center}
for $ \BV_1, \BV_2 $. That is, we let
\begin{align*}
\iota(\bra 0 \ket) &= \BV_1 \boxtimes \BV_2 \\
\iota(\bra 1 \ket) &= \BV_1 \boxtimes \BA \boxtimes \BV_2 \\
\iota(\bra 2 \ket) &= \BV_1 \boxtimes \BA \boxtimes \BA \boxtimes \BV_2,  
\end{align*}
and consider the $ 1 $-morphisms $ d_i = \iota(\partial_i) $ given by 
\begin{align*}
d_0 &= \rho \boxtimes \id \\
d_1 &= \id \boxtimes \lambda
\end{align*}
in degree $ 1 $, 
\begin{align*}
d_0 &= \rho \boxtimes \id \boxtimes \id \\
d_1 &= \id \boxtimes \otimes \boxtimes \id \\
d_2 &= \id \boxtimes \id \boxtimes \lambda
\end{align*}
in degree $ 2 $, and we set 
$$
\iota(\partial_{j - 1} \circ \partial_i) = \iota(\partial_i \circ \partial_j) = \iota(\partial_i) \circ \iota(\partial_j). 
$$
for $ i < j $. We also define 
$ \iota_{\bra 2 \ket, \bra 1 \ket, \bra 0 \ket}: \iota(\partial_i) \circ \iota(\partial_j) \Rightarrow \iota(\partial_i \circ \partial_j) $ 
to be the identity transformation for $ i < j $, and 
let $ \iota_{\bra 2 \ket, \bra 1 \ket, \bra 0 \ket}: \iota(\partial_{j - 1}) \circ \iota(\partial_i) \Rightarrow \iota(\partial_{j - 1} \circ \partial_i) $ 
be given by the associativity constraints, that is, 
\begin{align*} 
\alpha \boxtimes \id: d_0 \circ d_0 &\Rightarrow d_0 \circ d_1 \\
\id: d_1 \circ d_0 &\Rightarrow d_0 \circ d_2 \\
\id \boxtimes \alpha: d_1 \circ d_1 &\Rightarrow d_1 \circ d_2, 
\end{align*}
respectively. Together with $ \iota(1_{\bra j \ket}) = 1_{\iota(\bra j \ket)} $ and $ \iota_{\bra j \ket} = \id $ 
for $ j = 0,1,2 $ this determines $ \iota $ uniquely. 

According to Theorem \ref{bicolimitexistence} there exists a bicolimit $ \BV_1 \boxtimes_\BA \BV_2 \in C^* \Lin $ for this diagram, together with a 
transformation $ T: \iota \rightarrow \Delta(\BV_1 \boxtimes_\BA \BV_2) $ from $ \iota $ to the constant diagram with value $ \BV_1 \boxtimes_\BA \BV_2 $. 
In particular, there exists a nondegenerate 
linear $ * $-functor $ \boxtimes_\BA = T_{\bra 0 \ket}: \iota(\bra 0 \ket) = \BV_1 \boxtimes \BV_2 \rightarrow M(\BV_1 \boxtimes_\BA \BV_2) $ 
and a multiplier natural isomorphism 
$ \beta = T_{\bra 1 \ket,\bra 0 \ket}(\partial_1) \circ T_{\bra 1 \ket,\bra 0 \ket}(\partial_0)^{-1}: \boxtimes_\BA \circ (\rho \boxtimes \id) 
\Rightarrow \boxtimes_\BA \circ (\id \boxtimes \lambda) $. 

The following argument, verifying that this data defines a balanced tensor product, is folklore, but we shall carry out the details for the sake 
of completeness. In order to improve legibility we will abbreviate $ j = \bra j \ket $ in the sequel. 

It suffices to show that the category of transformations $ \iota \rightarrow \Delta(\BX) $ is equivalent 
to the category of $ \BA $-balanced functors $ \BV_1 \boxtimes \BV_2 \rightarrow M\BX $ for $ \BX \in C^* \Lin $ in a natural way. 
For this we shall define $ F: [I, C^* \Lin](\iota, \Delta(\BX)) \rightarrow C^* \Bilin_\BA(\BV_1, \BV_2; \BX) $ 
on objects by sending a transformation $ \Bf: \iota \rightarrow \Delta(\BX) $ to $ F(\Bf) = \Bf_0 $, together with the multiplier natural 
isomorphism $ \beta_{F(\Bf)} = \Bf_{1,0}(\partial_1) \circ \Bf_{1,0}(\partial_0)^{-1} $. 
On the level of morphisms we define $ F(\phi) = \phi_0 $.  

To verify that $ F(\Bf) $ and $ \beta_{F(\Bf)} $ yield an $ \BA $-balanced functor we need to check the defining relation in 
Definition \ref{defbalancedfunctor}. Using that $ \Bf $ is a transformation we calculate 
\begin{align*}
&(\beta_{F(\Bf)} * \id_{\id \boxtimes \id \boxtimes \lambda}) \circ (\beta_{F(\Bf)} * \id_{\rho \boxtimes \id \boxtimes \id}) \\
&= ((\Bf_{1,0}(\partial_1) \circ \Bf_{1,0}(\partial_0)^{-1}) * \id_{d_2}) \circ 
((\Bf_{1,0}(\partial_1) \circ \Bf_{1,0}(\partial_0)^{-1}) * \id_{d_0}) \\ 
&= (\Bf_{1,0}(\partial_1) * \id_{d_2}) \circ \Bf_{2,1}(\partial_2) \circ \Bf_{2,1}(\partial_0)^{-1} \circ (\Bf_{1,0}(\partial_0)^{-1} * \id_{d_0}) \\
&= \Bf_{2, 0}(\partial_1 \circ \partial_2) \circ \Bf_{2, 0}(\partial_0 \circ \partial_0)^{-1} \circ (\id_{\Bf_0} * \iota_{2,1,0}) \\
&= \Bf_{2, 0}(\partial_1 \circ \partial_1) \circ \Bf_{2, 0}(\partial_0 \circ \partial_1)^{-1} \circ (\id_{\Bf_0} * \iota_{2,1,0}) \\
&= (\id_{\Bf_0} * \iota_{2,1,0}) \circ (\Bf_{1,0}(\partial_1) * \id_{d_1}) \circ \Bf_{2,1}(\partial_1) \circ \Bf_{2, 0}(\partial_0 \circ \partial_1)^{-1} 
\circ (\id_{\Bf_0} * \iota_{2,1,0}) \\
&= (\id_{\Bf_0} * \iota_{2,1,0}) \circ (\Bf_{1,0}(\partial_1) * \id_{d_1}) \circ (\Bf_{1,0}(\partial_0)^{-1} * \id_{d_1}) 
\circ (\id_{\Bf_0} * \iota_{2,1,0}) \\
&= (\id_{F(\Bf)} * (\id \boxtimes \alpha)) \circ (\beta_{F(\Bf)} * \id_{\id \boxtimes \otimes \boxtimes \id}) \circ (\id_{F(\Bf)} * (\alpha \boxtimes \id))
\end{align*}
as required, showing that $ F $ is well-defined on objects. If $ \phi: \Bf \rightarrow \Bg $ 
is a modification then the relations $ (\phi_0 * \id_{d_j}) \circ \Bf_{1,0}(\partial_j) = \Bg_{1,0}(\partial_j) \circ \phi_1 $ for $ j = 0,1 $ imply
that $ F(\phi): F(\Bf) \Rightarrow F(\Bg) $ is an $ \BA $-balanced natural transformation. 
It follows that $ F $ defines a functor as stated. 

We claim that $ F $ is an equivalence of categories. 
Let $ \Bg: \BV_1 \boxtimes \BV_2 \rightarrow M\BX $ be an $ \BA $-balanced functor. We define a transformation $ \Bf: \iota \rightarrow \Delta(\BX) $ by setting 
\begin{align*}
\Bf_0 &= \Bg \\
\Bf_1 &= \Bg \circ d_0 \\
\Bf_2 &= \Bg \circ d_0 \circ d_0
\end{align*}
on objects. Moreover let 
\begin{align*} 
\Bf_{1,0}(\partial_0) &= \id: \Bf_1 \rightarrow \Bf_0 \circ d_0 \\
\Bf_{1,0}(\partial_1) &= \beta: \Bf_1 \rightarrow \Bf_0 \circ d_1 
\end{align*}
in degree $ 1 $, 
\begin{align*}
\Bf_{2,1}(\partial_0) &= \id: \Bf_2 = \Bg \circ d_0 \circ d_0 \rightarrow \Bg \circ d_0 \circ d_0 = \Bf_1 \circ d_0 \\
\Bf_{2,1}(\partial_1) &= \id_\Bg * (\alpha \boxtimes \id): \Bf_2 = \Bg \circ d_0 \circ d_0 \rightarrow \Bg \circ d_0 \circ d_1 = \Bf_1 \circ d_1 \\
\Bf_{2,1}(\partial_2) &= \beta * \id_{\rho \boxtimes \id \boxtimes \id}: \Bf_2 = \Bg \circ d_0 \circ d_0 \rightarrow \Bg \circ d_0 \circ d_2 = \Bf_1 \circ d_2 
\end{align*}
in degree $ 2 $ and set $ \Bf_{j,j}(1_j) = \id $ for $ j = 0,1,2 $. 
We also need to fix $ \Bf_{2,0} $ on the morphisms from level $ 2 $ to level $ 0 $ in $ I $,  
and in accordance with our choices for $ \iota $ we let 
\begin{align*}
\Bf_{2,0}(\partial_0 \circ \partial_0) = \Bf_{2,0}(\partial_0 \circ \partial_1) &= (\Bf_{1,0}(\partial_0) * \id_{d_1}) \circ \Bf_{2,1}(\partial_1) \\
\Bf_{2,0}(\partial_1 \circ \partial_0) = \Bf_{2,0}(\partial_0 \circ \partial_2) &= (\Bf_{1,0}(\partial_0) * \id_{d_2}) \circ \Bf_{2,1}(\partial_2) \\
\Bf_{2,0}(\partial_1 \circ \partial_1) = \Bf_{2,0}(\partial_1 \circ \partial_2) &= (\Bf_{1,0}(\partial_1) * \id_{d_2}) \circ \Bf_{2,1}(\partial_2).  
\end{align*}
Then the unit conditions for a transformation are trivially satisfied, and the remaining conditions read 
$$ 
\Bf_{2, 0}(\partial_i \circ \partial_j) = (\id_{\Bf_0} * \iota_{2,1,0}) \circ (\Bf_{1, 0}(\partial_i) * \id_{d_j}) \circ \Bf_{2, 1}(\partial_j)
$$
for all $ i, j $. For $ i < j $ these equalities hold by construction. In addition we have 
\begin{align*}
\Bf_{2, 0}(\partial_0 \circ \partial_0) &= \Bf_{2, 0}(\partial_0 \circ \partial_1) \\
&= (\Bf_{1, 0}(\partial_0) * \id_{d_1}) \circ \Bf_{2, 1}(\partial_1) \\
&= (\id * \id_{d_1}) \circ (\id_\Bg * (\alpha \boxtimes \id)) \\
&= \id_\Bg * (\alpha \boxtimes \id) \\
&= (\id_{\Bg} * (\alpha \boxtimes \id)) \circ (\id * \id_{d_0}) \circ \id \\
&= (\id_{\Bf_0} * \iota_{2,1,0}) \circ (\Bf_{1, 0}(\partial_0) * \id_{d_0}) \circ \Bf_{2, 1}(\partial_0) 
\end{align*}
as required, and similarly 
\begin{align*}
\Bf_{2, 0}(\partial_1 \circ \partial_0) &= \Bf_{2, 0}(\partial_0 \circ \partial_2) \\
&= (\Bf_{1,0}(\partial_0) * \id_{d_2}) \circ \Bf_{2,1}(\partial_2) \\
&= (\id * \id_{d_2}) \circ (\beta * \id_{\rho \boxtimes \id \boxtimes \id}) \\
&= \beta * \id_{\rho \boxtimes \id \boxtimes \id} \\
&= (\beta * \id_{d_0}) \circ \id \\
&= (\Bf_{1,0}(\partial_1) * \id_{d_0}) \circ \Bf_{2,1}(\partial_0). 
\end{align*}
Finally, the equality
\begin{align*}
\Bf_{2, 0}(\partial_1 \circ \partial_1) &= \Bf_{2, 0}(\partial_1 \circ \partial_2) \\
&= (\Bf_{1,0}(\partial_1) * \id_{d_2}) \circ \Bf_{2,1}(\partial_2) \\
&= (\beta * \id_{d_2}) \circ (\beta * \id_{\rho \boxtimes \id \boxtimes \id}) \\
&= (\id_\Bg * (\id \boxtimes \alpha)) \circ (\beta * \id_{\id \boxtimes \otimes \boxtimes \id}) \circ (\id_\Bg * (\alpha \boxtimes \id)) \\
&= (\id_{\Bf_0} * \iota_{2,1,0}) \circ (\Bf_{1,0}(\partial_1) * \id_{d_1}) \circ \Bf_{2,1}(\partial_1) 
\end{align*}
follows from the defining relation of the $ \BA $-balanced functor $ \Bg $. Hence $ \Bf $ is indeed a transformation,  
and since $ F(\Bf) = \Bg $ we conclude that $ F $ is essentially surjective. 

If $ \Bf, \Bg: \iota \rightarrow \Delta(\BX) $ are transformations and $ \phi: \Bf \rightarrow \Bg $ is a modification such that $ F(\phi) = \phi_0 = 0 $, then 
the modification property implies $ \phi_i = 0 $ for $ i = 0,1,2 $. 
It follows that $ F $ is faithful. Conversely, assume that $ \psi: F(\Bf) \rightarrow F(\Bg) $ is an $ \BA $-balanced transformation. 
Define $ \phi: \Bf \rightarrow \Bg $ by 
\begin{align*} 
\phi_0 &= \psi \\ 
\phi_1 &= \Bg_{1,0}(\partial_0)^{-1} \circ (\psi * \id_{d_0}) \circ \Bf_{1,0}(\partial_0) \\
\phi_2 &= \Bg_{2,1}(\partial_0)^{-1} \circ (\Bg_{1,0}(\partial_0)^{-1} * \id_{d_0}) 
\circ (\psi * \id_{d_0 \circ d_0}) \circ (\Bf_{1,0}(\partial_0) * \id_{d_0}) \circ \Bf_{2,1}(\partial_0).  
\end{align*}
Then $ (\phi_0 * \id_{d_j}) \circ \Bf_{1,0}(\partial_j) = \Bg_{1,0}(\partial_j) \circ \phi_1 $ for $ j = 0 $ by definition, 
and for $ j = 1 $ this relation follows form the fact that $ \psi $ is $ \BA $-balanced, that is, 
\begin{align*}
\Bg_{1,0}(\partial_1)^{-1} \circ (\psi * \id_{d_1}) \circ \Bf_{1,0}(\partial_1) 
&= \Bg_{1,0}(\partial_1)^{-1} \circ (\psi * \id_{d_1}) \circ \beta_{F(\Bf)} \circ \Bf_{1,0}(\partial_0) \\
&= \Bg_{1,0}(\partial_1)^{-1} \circ \beta_{F(\Bg)} \circ (\psi * \id_{d_0}) \circ \Bf_{1,0}(\partial_0) \\
&= \Bg_{1,0}(\partial_0)^{-1} \circ (\psi * \id_{d_0}) \circ \Bf_{1,0}(\partial_0). 
\end{align*}
Similarly, for $ j = 0 $ we obtain the relation $ (\phi_1 * \id_{d_j}) \circ \Bf_{2,1}(\partial_j) = \Bg_{2,1}(\partial_j) \circ \phi_2 $ 
directly from the definition. To check the case $ j = 1 $ note that 
\begin{align*}
\phi_2 &= \Bg_{2,0}(\partial_0 \circ \partial_0)^{-1} \circ (\id_{\Bg_0} * \iota_{2,1,0})  
\circ (\psi * \id_{d_0 \circ d_0}) \circ (\id_{\Bf_0} * \iota_{2,1,0}^{-1}) \circ \Bf_{2,0}(\partial_0 \circ \partial_0) \\
&= \Bg_{2,0}(\partial_0 \circ \partial_1)^{-1} 
\circ (\psi * \id_{d_0 \circ d_1}) \circ \Bf_{2,0}(\partial_0 \circ \partial_1) \\
&= \Bg_{2,1}(\partial_1)^{-1} \circ (\Bg_{1,0}(\partial_0)^{-1} * \id_{d_1}) 
\circ (\psi * \id_{d_0 \circ d_1}) \circ (\Bf_{1,0}(\partial_0) * \id_{d_1}) \circ \Bf_{2,1}(\partial_1).  
\end{align*} 
For $ j = 2 $ we calculate
\begin{align*}
\phi_2 &= \Bg_{2,1}(\partial_0)^{-1} \circ (\Bg_{1,0}(\partial_0)^{-1} * \id_{d_0}) 
\circ (\psi * \id_{d_0 \circ d_0}) \circ (\Bf_{1,0}(\partial_0) * \id_{d_0}) \circ \Bf_{2,1}(\partial_0) \\
&= \Bg_{2,1}(\partial_0)^{-1} \circ (\Bg_{1,0}(\partial_1)^{-1} * \id_{d_0}) 
\circ (\psi * \id_{d_1 \circ d_0}) \circ (\Bf_{1,0}(\partial_1) * \id_{d_0}) \circ \Bf_{2,1}(\partial_0) \\
&= \Bg_{2,0}(\partial_1 \circ \partial_0)^{-1}  
\circ (\psi * \id_{d_1 \circ d_0}) \circ \Bf_{2,0}(\partial_1 \circ \partial_0) \\
&= \Bg_{2,0}(\partial_0 \circ \partial_2)^{-1} 
\circ (\psi * \id_{d_0 \circ d_2}) \circ \Bf_{2,0}(\partial_0 \circ \partial_2) \\
&= \Bg_{2,1}(\partial_2)^{-1} \circ (\Bg_{1,0}(\partial_0)^{-1} * \id_{d_2}) 
\circ (\psi * \id_{d_0 \circ d_2}) \circ (\Bf_{1,0}(\partial_0) * \id_{d_2}) \circ \Bf_{2,1}(\partial_2).  
\end{align*} 
It follows that $ \phi: \Bf \rightarrow \Bg $ is a modification, and we have $ F(\phi) = \psi $ by construction. This shows that $ F $ is full, 
and finishes the proof. \qed 

Let us conclude our discussion by describing two simple examples of balanced tensor products. 

\begin{lemma} 
Let $ \BA $ be a countably additive $ C^* $-tensor category. If $ \BV $ is a left $ \BA $-module category then the module category structure induces an 
equivalence of $ \BA $-module categories
$$ 
\BA \boxtimes_\BA \BV \simeq \BV. 
$$
An analogous statement holds for right module categories. 
\end{lemma} 

\proof The nondegenerate linear $ * $-functor $ \lambda: \BA \boxtimes_{\max} \BV \rightarrow \BV $ giving the module category structure is 
$ \BA $-balanced, so that it induces a nondegenerate linear $ * $-functor $ \lambda: \BA \boxtimes_\BA \BV \rightarrow M\BV $. We get a 
linear $ * $-functor $ \eta: \BV \rightarrow M(\BA \boxtimes_\BA \BV) $ 
induced by the functor $ \BV \rightarrow M(\BA \boxtimes_{\max} \BV) $ given by $ \eta(V) = 1 \boxtimes V $. 
The composition $ \lambda \circ \eta $ is easily seen to be naturally isomorphic to the identity. 
In the opposite direction we get $ (\eta \circ \lambda)(X \boxtimes V) = 1 \boxtimes (X \otimes V) \cong X \boxtimes V $ 
for $ X \in \BA $ and $ V \in \BV $ by balancedness, and it follows that $ \eta \circ \lambda $ is naturally isomorphic to the identity. \qed 

Recall the infinite matrix $ C^* $-tensor category over a countably additive $ C^* $-tensor category and its module categories described 
in Example \ref{moduleexamples} c).  

\begin{prop} 
Let $ \BA $ be a countably additive $ C^* $-tensor category. Moreover let $ \KH_\BA $ be the infinite matrix $ C^* $-tensor category over $ \BA $ 
and let $ \HH_\BA $ and $ \HH_\BA^\vee $ be the associated column and row bimodule categories, respectively. Then we have equivalences 
$$
\HH_\BA^\vee \boxtimes_{\KH_\BA} \HH_\BA \simeq \BA, \qquad \HH_\BA \boxtimes_\BA \HH_\BA^\vee \simeq \KH_\BA 
$$
of $ \BA $-$ \BA $ and $ \KH_\BA $-$ \KH_\BA $-bimodule categories, respectively. 
\end{prop} 

\proof We obtain balanced functors $ \HH_\BA^\vee \times \HH_\BA \rightarrow \BA $ and $ \HH_\BA \boxtimes \HH_\BA^\vee \rightarrow \KH_\BA $ 
using categorical matrix multiplication. These are naturally functors of $ \BA $-$ \BA $-bimodule categories and $ \KH_\BA $-$ \KH_\BA $-bimodule
categories, respectively. It is straightforward to check that both satisfy the defining property of a balanced tensor product. \qed

\appendix

\section{Bicategories} \label{appendix} 

We collect some definitions and facts related to bicategories. For more information see for instance \cite{Benabou}, \cite{Fiorebook}. 
In the sequel all categories are unital. 

\subsection*{Bicategories, homomorphisms, transformations, modifications} 

Let us start by recalling the definition of a bicategory.  

\begin{definition} \label{defbicat}
A bicategory $ \B $ consists of
\begin{bnum}
\item[$ \bullet $] a class $ \B_0 $ of objects or 0-cells. 
\item[$ \bullet $] categories $ \B(\BX, \BY) $ for all objects $ \BX, \BY $ in $ \B_0 $. The objects of these 
categories are called the 1-cells or 1-morphisms of $ \B $, and we write $ \B_1 $ for the collection of all $ 1 $-morphisms. 
We shall typically write $ 1 $-morphisms using letters $ \Bf, \Bg, \Bh $, or $ \Bf: \BX \rightarrow \BY $ if the source and target objects 
shall be indicated. Morphisms in $ \B(\BX, \BY) $ are called 2-cells or 2-morphisms, and we write $ \B_2 $ for the collection of 
all of them. 
Typical $ 2 $-morphisms will be denoted $ \alpha, \beta, \gamma $. We will also write $ \alpha: \Bf \rightarrow \Bg $ or $ \alpha: \Bf \Rightarrow \Bg $ 
if the source and target $ 1 $-cells shall be indicated. The composition of $ 2 $-morphisms $ \alpha: \Bf \rightarrow \Bg, \beta: \Bg \rightarrow \Bh $ 
is called vertical composition and denoted by $ \beta \circ \alpha: \Bf \rightarrow \Bh $. 
\item[$ \bullet $] composition functors 
$$ 
c_{\BX, \BY, \BZ}: \B(\BX, \BY) \times \B(\BY,\BZ) \rightarrow \B(\BX, \BZ) 
$$ 
for all objects $ \BX, \BY, \BZ \in \B_0 $, called horizontal composition. 
We shall write $ (\Bf, \Bg) \mapsto \Bg \circ \Bf = c_{\BX, \BY, \BZ}(\Bf, \Bg) $ as for the composition of maps between sets. We follow the standard 
convention to write $ \beta * \alpha = c_{\BX, \BY, \BZ}(\alpha, \beta): (\Bg \circ \Bf) \rightarrow (\Bk \circ \Bh) $ for the horizontal composition 
of $ 2 $-morphisms $ \alpha: \Bf \rightarrow \Bh, \beta: \Bg \rightarrow \Bk $. 
\item[$ \bullet $] an object $ 1_\BX \in \B(\BX, \BX) $ for each $ \BX \in \B_0 $ called the identity of $ \BX $, sometimes 
identified with the image of a functor $ \star \rightarrow \B(\BX, \BX) $ from the category with one object and one morphism. 
\item[$ \bullet $] natural isomorphisms 
$$ 
\alpha_{\BW, \BX, \BY, \BZ}: c_{\BW, \BX, \BZ} \circ (\id \times c_{\BX, \BY, \BZ}) \rightarrow c_{\BW, \BY, \BZ} \circ (c_{\BW, \BX, \BY} \times \id) 
$$ 
for all $ \BW, \BX, \BY, \BZ \in \B_0 $, written 
$$ 
\alpha_{\BW, \BX, \BY, \BZ}: (\Bh \circ \Bg) \circ \Bf \rightarrow \Bh \circ (\Bg \circ \Bf) 
$$ 
on the level of objects, called associators. 
\item[$ \bullet $] natural isomorphisms $ \rho_{\BX, \BY}: c_{\BX, \BX, \BY} \circ (1_\BX \times \id) \rightarrow \id $ 
and $ \lambda_{\BX, \BY}: c_{\BX, \BY, \BY} \circ (\id \times 1_\BY) \rightarrow \id $ in $ \B(\BX, \BY) $ for all $ \BX, \BY \in \B_0 $, 
called left and right unitors, and written $ \rho_{\BX, \BY}(\Bf): \Bf \circ 1_\BX \rightarrow \Bf $ 
and $ \lambda_{\BX, \BY}(\Bf): 1_\BY \circ \Bf \rightarrow \Bf $ for $ \Bf \in \B(\BX, \BY) $, respectively. 
\end{bnum}
These data are supposed to satisfy the following conditions. 
\begin{bnum}
\item[$\bullet$] (Associativity constraints) For all 
1-morphisms $ \Bf: \BV \rightarrow \BW, \Bg: \BW \rightarrow \BX, \Bh: \BX \rightarrow \BY, \Bk: \BY \rightarrow \BZ $ 
the diagram
\begin{center}
\begin{tikzcd}[column sep = -5mm]
			&&
			(\Bk \circ \Bh) \circ (\Bg \circ \Bf)
			\arrow[drr,"\alpha_{\BV, \BX, \BY, \BZ}"]
			&&
			\\
			((\Bk \circ \Bh) \circ \Bg) \circ \Bf
			\arrow[urr,"\alpha_{\BV, \BW, \BX, \BZ}"]
			\arrow[dr, swap, "\alpha_{\BW, \BX, \BY, \BZ} * \id_\Bf "]
			&&&&
			\Bk \circ (\Bh \circ (\Bg \circ \Bf))
			\\
			&
			(\Bk \circ (\Bh \circ \Bg)) \circ \Bf
			\arrow[rr, swap, "\alpha_{\BV, \BW, \BY, \BZ}"]
			&&
			\Bk \circ ((\Bh \circ \Bg) \circ \Bf)
			\arrow[ur, swap, "\id_\Bk * \alpha_{\BV, \BW, \BX, \BY}"]
			&
\end{tikzcd}
\end{center}
is commutative. 
\item[$\bullet$] (Unit constraints) For all 1-morphisms $ \Bf: \BX \rightarrow \BY, \Bg: \BY \rightarrow \BZ $ the diagram 
\begin{center}
\begin{tikzcd}[column sep = 15mm]
			(\Bg \circ 1_\BY) \circ \Bf
			\arrow[rr, "\alpha_{\BX, \BY, \BY, \BZ}"]
			\arrow[dr, swap, "\rho_{\BY, \BZ}(\Bg) * \id_\Bf"]
			&&
			\Bg \circ (1_\BY \circ \Bf)
			\arrow[dl, "\id_\Bg * \lambda_{\BX, \BY}(\Bf)"]
			\\
			&
			\Bg \circ \Bf
\end{tikzcd}
\end{center}
is commutative.  
\end{bnum} 
\end{definition}

If all associators $ \alpha $ and all unitors $ \lambda, \rho $ in the bicategory $ \B $ are identity transformations then $ \B $ is called strict, 
or a $ 2 $-category. Any bicategory is equivalent to a $ 2 $-category, for the notion of equivalence to be explained below. 

To improve legibility we will sometimes omit subscripts in the natural transformations $ \alpha, \lambda, \rho $, and write for 
instance $ \alpha $ instead of $ \alpha_{\BW, \BX, \BY, \BZ} $. 

\begin{definition} \label{defbicathom}
Let $ \B, \C $ be bicategories. A homomorphism $ f: \B \rightarrow \C $ consists of 
\begin{bnum}
\item[$\bullet$] a map $ f: \B_0 \rightarrow \C_0 $ between the objects of $ \B $ and $ \C $, 
\item[$\bullet$] functors $ f_{\BX, \BY}: \B(\BX, \BY) \rightarrow \C(f(\BX), f(\BY)) $ for all $ \BX, \BY \in \B $, 
\item[$\bullet$] natural isomorphisms $ f_{\BX, \BY, \BZ}: c_{f(\BX), f(\BY), f(\BZ)} \circ (f_{\BY, \BZ} \times f_{\BX, \BY}) \rightarrow 
f_{\BX, \BZ} \circ c_{\BX, \BY, \BZ} $ and $ f_\BX: 1_{f(\BX)} \rightarrow f_{\BX, \BX}(1_\BX) $
\end{bnum} 
such that the following conditions hold. 
\begin{bnum} 
\item[$\bullet$] For all $ \Bf: \BW \rightarrow \BX, \Bg: \BX \rightarrow \BY, \Bh: \BY \rightarrow \BZ $ the diagram
\begin{center}
\begin{tikzcd}[column sep = 0mm]
		  &
			(f(\Bh) \circ f(\Bg)) \circ f(\Bf)
			\arrow[dl, swap, "\alpha"]
			\arrow[dr,"f_{\BX, \BY, \BZ} * \id_{f(\Bf)}"]
			&
			\\
			f(\Bh) \circ (f(\Bg) \circ f(\Bf))
			\arrow[d, swap, "\id_{f(\Bh)} * f_{\BW, \BY, \BZ}"]
			&&
			f(\Bh \circ \Bg) \circ f(\Bf)
			\arrow[d, "f_{\BW, \BX, \BZ}"]
			\\
			f(\Bh) \circ f(\Bg \circ \Bf)
			\arrow[dr, swap, "f_{\BW, \BY, \BZ}"]
      &&
			f((\Bh \circ \Bg) \circ \Bf)
			\arrow[dl, "f(\alpha)"]
			\\
			&
			f(\Bh \circ (\Bg \circ \Bf))
			&
\end{tikzcd}
\end{center}
is commutative. 
\item[$\bullet$] For all $ \Bf: \BX \rightarrow \BY $ the diagrams 
\begin{center}
\begin{tikzcd}[column sep = 5mm]
      &
			f(\Bf) \circ 1_{f(\BX)}
    	\arrow[dl, swap, "\id_{f(\Bf)} * f_\BX"]
			\arrow[dr, "\rho"]
			&&
			1_{f(\BY)} \circ f(\Bf) 
			\arrow[dl, swap, "\lambda"]
			\arrow[dr, "f_\BY * \id_{f(\Bf)}"]
			&
			\\
			f(\Bf) \circ f(1_\BX)
			\arrow[dr, swap, "f_{\BX, \BX, \BY}"]
			&
			&
			f(\Bf)
			&
			&
			f(1_\BY) \circ f(\Bf)
			\arrow[dl, "f_{\BX, \BY, \BY}"]
			\\
			&
			f(\Bf \circ 1_\BX)
			\arrow[ur, swap, "f(\rho)"]
			&&
			f(1_\BY \circ \Bf)
      \arrow[ul, "f(\lambda)"]
			&
\end{tikzcd}
\end{center}
are commutative. 
\end{bnum}
\end{definition}

Homomorphisms of bicategories in the sense of Definition \ref{defbicathom} are also called weak $ 2 $-functors or pseudofunctors. A 
homomorphism $ f: \B \rightarrow \C $ is called strict if $ f_{\BX, \BY, \BZ} $ and $ f_\BX $ are the identity for all $ \BX, \BY, \BZ \in \B $. 

\begin{definition} \label{defbicattrans}
Let $ \B, \C $ be bicategories and $ f,g: \B \rightarrow \C $ homomorphisms. A transformation $ \sigma: f \rightarrow g $ consists of 
\begin{bnum}
\item[$\bullet$] a 1-morphism $ \sigma_\BX: f(\BX) \rightarrow g(\BX) $ for every object $ \BX $ of $ \B $, 
\item[$\bullet$] natural isomorphisms $ \sigma_{\BX, \BY}: (\sigma_\BX)^* g_{\BX, \BY} \rightarrow (\sigma_\BY)_* f_{\BX, \BY} $ 
for all objects $ \BX, \BY $ of $ \B $, so that we have invertible 
2- morphisms 
$$ 
\sigma_{\BX, \BY}(\Bf): g(\Bf) \circ \sigma_\BX \rightarrow \sigma_\BY \circ f(\Bf) 
$$ 
for all $ \Bf \in \B(\BX, \BY) $. 
\end{bnum} 
These data are required to satisfy the following conditions. 
\begin{bnum} 
\item[$\bullet$] For all $ \Bf: \BX \rightarrow \BY, \Bg: \BY \rightarrow \BZ $ the diagram
\begin{center}
\begin{tikzcd}[column sep = 0mm]
		  &
			g(\Bg \circ \Bf) \circ \sigma_\BX
			\arrow[dr,"\sigma_{\BX, \BZ}(\Bg \circ \Bf)"]
			&
			\\
			(g(\Bg) \circ g(\Bf)) \circ \sigma_\BX
			\arrow[ur, "g_{\BX, \BY, \BZ} * \id_{\sigma_\BX}"]
			\arrow[d, swap, "\alpha_{f(\BX), g(\BX), g(\BY), g(\BZ)}"]
			&&
			\sigma_\BZ \circ f(\Bg \circ \Bf)
			\\
			g(\Bg) \circ (g(\Bf) \circ \sigma_\BX)
			\arrow[d, swap, "\id_{g(\Bg)} * \sigma_{\BX, \BY}(\Bf)"]
			&&
			\sigma_\BZ \circ (f(\Bg) \circ f(\Bf))
			\arrow[u, swap, "\id_{\sigma_\BZ} * f_{\BX, \BY, \BZ}"]
			\\
			g(\Bg) \circ (\sigma_\BY \circ f(\Bf))
			\arrow[dr, swap, "\alpha_{\BX, \BY, \BY, \BZ}^{-1}"]
      &&
			(\sigma_\BZ \circ f(\Bg)) \circ f(\Bf)
			\arrow[u, swap, "\alpha_{f(\BX), f(\BY), f(\BZ), g(\BZ)}"]
			\\
			&
			(g(\Bg) \circ \sigma_\BY) \circ f(\Bf)
			\arrow[ur, swap, "\sigma_{\BY, \BZ}(\Bg) * \id_{f(\Bf)}"]
			&
\end{tikzcd}
\end{center}
is commutative. 
\item[$\bullet$] For all objects $ \BX $ in $ \B $ the diagram
\begin{center}
\begin{tikzcd}[column sep = -7mm]
			&&
			\qquad \sigma_\BX \qquad 
			\arrow[drr,"\rho_{g(\BX), f(\BX)}^{-1}", start anchor = {[xshift = -3mm, yshift = 1mm]}, end anchor = {[xshift = 2mm]}]
			&&
			\\
			1_{g(\BX)} \circ \sigma_\BX
			\arrow[urr,"\lambda_{g(\BX), f(\BX)}", start anchor = {[xshift = -2mm]}, end anchor = {[xshift = 3mm, yshift = 1mm]}]
			\arrow[dr, swap, "g_\BX * \id_{\sigma_\BX}"]
			&&&&
			\sigma_\BX \circ 1_{f(\BX)}
			\arrow[dl, "\id_{\sigma_\BX} * f_\BX"]
			\\
			&
			g(1_\BX) \circ \sigma_\BX
			\arrow[rr, swap, "\sigma_{1_\BX}"]
			&&
			\sigma_\BX \circ f(1_\BX) 
			&
\end{tikzcd}
\end{center}
is commutative. 
\end{bnum} 
\end{definition}

Transformations are sometimes also referred to as pseudonatural transformations. 

\begin{definition} \label{defbicatmod}
Let $ \B, \C $ be bicategories, $ f,g: \B \rightarrow \C $ be homomorphisms and $ \sigma, \tau: f \rightarrow g $ be transformations. 
A modification $ \Gamma: \sigma \rightarrow \tau $ consists of $ 2 $-morphisms $ \Gamma_\BX: \sigma_\BX \rightarrow \tau_\BX $ for all objects $ \BX \in \B $ 
such that the diagram 
\begin{center}
\begin{tikzcd}[column sep = 15mm]
			g(\Bf) \circ \sigma_\BX
			\arrow[r, "\id_{g(\Bf)} * \Gamma_\BX"]
			\arrow[d, swap, "\sigma_{\BX, \BY}(\Bf)"]
			&
			g(\Bf) \circ \tau_\BX
			\arrow[d, "\tau_{\BX, \BY}(\Bf)"]
			\\
			\sigma_\BY \circ f(\Bf)
			\arrow[r, swap, "\Gamma_\BY * \id_{f(\Bf)}"]
			&
			\tau_\BY \circ f(\Bf)
\end{tikzcd}
\end{center}
is commutative for every 1-morphism $ \Bf: \BX \rightarrow \BY $. 
\end{definition} 

There is a bicategory $ \Bicat(\B, \C) $ of homomorphism for bicategories $ \B $ and $ \C $, with objects 
homomorphisms from $ \B $ to $ \C $, transformations as 1-morphisms, and modifications as 2-morphisms. 

Two objects $ \BX, \BY $ in a bicategory $ \B $ are called equivalent if there exist $ 1 $-morphisms $ \Bf: \BX \rightarrow \BY, \Bg: \BY \rightarrow \BX $ 
such that their mutual compositions are isomorphic to the identities. 

Let us also discuss the composition of homomorphisms between bicategories. 

\begin{definition} 
Let $ f: \B \rightarrow \C $ and $ g: \C \rightarrow \D $ be homomorphims of bicategories. Then their composition $ g \otimes f: \B \rightarrow \D $ is defined by 
\begin{bnum}
\item[$\bullet$] the map $ g \otimes f: \B \rightarrow \D $ on objects given by the ordinary composition 
of $ f: \B_0 \rightarrow \C_0 $ and $ g: \C_0 \rightarrow \D_0 $,  
\item[$\bullet$] the functors $ (g \otimes f)_{\BX, \BY}: \B(\BX, \BY) \rightarrow \D((g \otimes f)(\BX), (g \otimes f)(\BY)) $ 
obtained as composition of $ f_{\BX, \BY}: \B(\BX, \BY) \rightarrow \C(f(\BX), f(\BY)) $ 
and $ g_{f(\BX), f(\BY)}: \C(f(\BX), f(\BY)) \rightarrow \D((g \otimes f)(\BX), (g \otimes f)(\BY)) $ for all $ \BX, \BY \in \B $, 
\item[$\bullet$] the natural isomorphisms $ (g \otimes f)_{\BX, \BY, \BZ}: c_{(g \otimes f)(\BX), (g \otimes f)(\BY), (g \otimes f)(\BZ)} \otimes 
((g \otimes f)_{\BY, \BZ} \times (g \otimes f)_{\BX, \BY}) \rightarrow (g \otimes f)_{\BX, \BZ} \otimes c_{\BX, \BY, \BZ} $ 
given by 
\begin{center}
\begin{tikzcd}[sep = huge]
		c_{(g \otimes f)(\BX), (g \otimes f)(\BY), (g \otimes f)(\BZ)} \otimes (g_{f(\BY), f(\BZ)} \times g_{f(\BX), f(\BY)}) \otimes (f_{\BY, \BZ} \times f_{\BX, \BY})
		\arrow[d, "g_{f(\BX), f(\BY), f(\BZ)} \otimes (f_{\BY, \BZ} \times f_{\BX, \BY})"]
		\\
  	g_{f(\BX), f(\BZ)} \otimes c_{f(\BX), f(\BY), f(\BZ)} \otimes (f_{\BY, \BZ} \times f_{\BX, \BY})
		\arrow[d, "g_{f(\BX), f(\BZ)} \otimes f_{\BX, \BY, \BZ}"]
		\\
		g_{f(\BX), f(\BZ)} \otimes f_{\BX, \BZ} \otimes c_{\BX, \BY, \BZ}
\end{tikzcd}
\end{center} 
and $ (g \otimes f)_\BX: 1_{g \otimes f(\BX)} \rightarrow (g \otimes f)_{\BX, \BX}(1_\BX) $ given by 
\begin{center}
\begin{tikzcd}[sep = huge]
		1_{g \otimes f(\BX)}
		\arrow[r, "g_{f(\BX)}"]
		&
  	g_{f(\BX), f(\BX)}(1_{f(\BX)})
		\arrow[r, "g_{f(\BX), f(\BX)}(f_\BX)"]
		&
		(g \otimes f)_{\BX, \BX}(1_\BX). 
\end{tikzcd}
\end{center} 
\end{bnum} 
\end{definition} 

The composition of two homomorphisms is again a homomorphism of bicategories, composition of homomorphisms is strictly associative, and the identity 
homomorphisms $ I_\B: \B \rightarrow \B $ are strict identities with respect to composition. 
In this way one obtains an ordinary category of bicategories and homomorphisms. 

\begin{definition}
Let $ \B, \C $ be bicategories. Then $ \B $ and $ \C $ are called biequivalent if there exist homomorphisms $ f: \B \rightarrow \C, g: \C \rightarrow \B $ 
and equivalences $ g \otimes f \simeq \id_\D $ in $ \Bicat(\B, \B) $ and $ f \otimes g \simeq \id_\C $ in $ \Bicat(\C, \C) $, respectively. 
\end{definition}

\subsection*{Biadjunctions} 

Biadjoint homomorphisms are the bicategorical analogue of adjoint functors. 

\begin{definition} 
Let $ \B, \C $ be bicategories and let $ f: \B \rightarrow \C, g: \C \rightarrow \B $ be homomorphisms. 
Then $ f $ is left biadjoint to $ g $, or equivalently $ g $ is right biadjoint to $ f $, if for all $ \BX \in \B, \BY \in \C $ there exist equivalences 
of categories 
$$
\C(f(\BX), \BY) \cong \B(\BX, g(\BY)),  
$$
pseudo-natural both in $ \BX $ and $ \BY $. 
\end{definition} 

Biadjunctions can also be described in terms of \emph{biuniversal arrows} in the following sense. 

\begin{definition} 
Let $ g: \C \rightarrow \B $ be a homomorphism of bicategories and let $ \BX \in \B $ be an object. 
An object $ f(\BX) \in \C $ together with a $ 1 $-morphism $ \eta_\BX \in \B(\BX, g(f(\BX))) $ is called a biuniversal arrow if 
for every object $ \BY \in \C $ the functor
$$ 
\sigma_\BY: \C(f(\BX), \BY) \rightarrow \B(\BX, g(\BY)) 
$$ 
defined by $ \sigma_\BY(\Bh) = g(\Bh) \circ \eta_\BX $ and $ \sigma_\BY(\gamma) = g(\gamma) * 1_{\eta_\BX} $ is an equivalence of categories. 
\end{definition} 

A family of biuniversal arrows determines a biadjunction as follows. 

\begin{theorem} \label{biadjointablechar}
Let $ g: \C \rightarrow \B $ be a homomorphism of bicategories. Then there exists a left biadjoint of $ g $ iff for every $ \BX \in \B $ 
there exists an object $ f(\BX) \in \C $ and a biuniversal arrow $ \eta_\BX: \BX \rightarrow g(f(\BX)) $. Moreover the left biadjoint 
homomorphism $ f: \B \rightarrow \C $ may be chosen to send $ \BX $ to $ f(\BX) $ on the level of objects. 
\end{theorem} 

For a proof of Theorem \ref{biadjointablechar} in the case of interest to us, namely if $ \B, \C $ are $ 2 $-categories, see chapter 9 in \cite{Fiorebook}. 

\subsection*{Bicategorical Yoneda Lemma} 

Let us discuss the analogue of the Yoneda Lemma in the setting of bicategories. 

If $ \B $ is a bicategory and $ \BX \in \B $ then there exists a homomorphism of bicategories $ Y_\BX: \B \rightarrow \Cat $ given by 
$$
Y_\BX(\BY) = \B(\BX, \BY) 
$$
on objects. For a $ 1 $-morphism $ \Bf: \BY \rightarrow \BZ $ one defines $ Y_\BX(\Bf): \B(\BX, \BY) \rightarrow \B(\BX, \BZ) $ 
by $ Y_\BX(\Bf)(\Bh) = \Bf \circ \Bh $, and on $ 2 $-morphisms one defines $ Y_\BX(\sigma) = \id_\Bf * \sigma $. 

We refer to $ Y_\BX $ as the representable homomorphism, or representable $ 2 $-presheaf, defined by $ \BX \in \B $. 
More generally, a homomorphism $ f: \B \rightarrow \Cat $ is called representable if there exists an object $ \BX \in \B $ and a 
transformation $ \sigma: f \rightarrow Y_\BX $  which is invertible up to modifications. 

Recall that $ \Bicat(\B, \C) $ is naturally a bicategory for all bicategories $ \B, \C $. If the target bicategory $ \C $ is a $ 2 $-category, 
then $ \Bicat(\B, \C) $ is a $ 2 $-category. 

\begin{theorem}[Bicategorical Yoneda] \label{bicategoricalyoneda}
Let $ \B $ be a bicategory. For every homomorphism $ f: \B \rightarrow \Cat $ and every $ \BX \in \B $ there is an equivalence of categories 
$$
\Bicat(\B, \Cat)(Y_\BX, f) \simeq f(\BX),  
$$
pseudonatural in $ f $ and $ \BX $. 
\end{theorem} 

For a detailed proof of Theorem \ref{bicategoricalyoneda} see \cite{Bakovicyoneda}.

\subsection*{Bilimits and bicolimits} 

Let us finally review the general definition of bicolimits in bicategories. The definition of bilimits is analogous. 

Let $ I $ be a small category. We will also view $ I $ as a $ 2 $-category with only identity $ 2 $-morphisms. 
If $ \B $ is a bicategory, then an $ I $-diagram in $ \B $ is a homomorphism $ \iota: I \rightarrow \B $. 
We will also use the notation $ (\BV_i)_{i \in I} $ for such an $ I $-diagram $ \iota $, denoting $ \BV_i = \iota(i) $ the assignment on the level of objects, 
and suppressing the remaining data. 

If $ \BW \in \B $ we obtain the constant $ I $-diagram $ \Delta(\BW) $ by setting $ \Delta(\BW)(i) = \BW $ and $ \Delta(\BW)(i \rightarrow j) = 1_\BW $ 
for all $ i \rightarrow j $ in $ I $. This defines a homomorphism of bicategories $ \Delta: \B \rightarrow \Bicat(I, \B) = [I, \B] $. 

\begin{definition} 
Let $ I $ be a small category and let $ (\BV_i)_{i \in I} $ be an $ I $-diagram in $ \B $. 
A bicolimit of $ (\BV_i)_{i \in I} $ is an object $ \varinjlim_{i \in I} \BV_i $ in $ \B $ together with an equivalence of categories
$$
\B(\varinjlim_{i \in I} \BV_i, \BX) \simeq [I, \B]((\BV_i)_{i \in I}, \Delta(\BX)),  
$$
pseudonatural in $ \BX $. 
\end{definition} 

It follows from the bicategorical Yoneda Lemma that a bicolimit is uniquely determined up to equivalence. 
If the bicolimits of all $ I $-diagrams in $ \B $ exist then $ \varinjlim_{i \in I} $ defines a left biadjoint to the constant $ I $-diagram 
homomorphism $ \Delta: \B \rightarrow \Bicat(I, \B) $.

\bibliographystyle{plain}

\bibliography{cvoigt}

\end{document}